\documentclass[a4paper,11pt]{article}
\usepackage{amsmath}
\usepackage{amssymb}
\usepackage{color}
\usepackage{graphicx}
\usepackage{epsfig}
\usepackage{epstopdf}
\usepackage[applemac]{inputenc}
\usepackage{url}
\usepackage{placeins}
\definecolor{labelkey}{rgb}{0.6,0,1}
\usepackage{epsfig}

\catcode`@=11
\@addtoreset{equation}{section}

\def\D{{\Delta}}
\def\1{1}

\newcommand{\R}{{\mathbb R}}

\newcommand{\G}{{\mathcal{G}}}

\newcommand{\qed}{ \begin{flushright}{\hfill\ifhmode\unskip\nobreak\fi
\ifmmode\ifinner\else\hskip 5pt\fi\fi
\hbox{\hskip 5pt\lower 1.5pt\hbox{\vrule width 0.2pt
\vbox{\hrule width 4pt height 0.2pt \vskip 7.2pt
\hrule width 4pt height 0.2pt}\unskip\vrule width
0.2pt}\hskip 0pt}} \end{flushright}}

\begin{document}
\title{A Semi-Lagrangian Scheme with Radial Basis Approximation for Surface Reconstruction \footnote{This research has been partially funded by Sapienza Universit\`a di Roma, Universit\`a Roma Tre and INdAM--GNCS}}
\author{E. Carlini \footnote{Dipartimento di Matematica ``G. Castelnuovo'', Sapienza Universit\`a di Roma, P.le A. Moro, 5, 00185 Roma (Italy), e-mail: {\tt carlini@mat.uniroma1.it}} \and R. Ferretti \footnote{Corresponding Author, Dipartimento di Matematica e Fisica, Universit\`a Roma Tre, L.go S. Leonardo Murialdo, 1, 00146 Roma (Italy), e-mail: {\tt ferretti@mat.uniroma3.it}, fax: +39-0657338080}}
\date{\today}
\maketitle
\begin{abstract} 
We propose a Semi-Lagrangian scheme coupled with Radial Basis Function interpolation for approximating a curvature-related level set model, which has been  proposed by Zhao et al. in \cite{ZOMK} to reconstruct unknown surfaces from sparse data sets. 
The main advantages of the proposed scheme are the possibility to solve the level set method on unstructured grids, as well as to concentrate the reconstruction points in the neighbourhood of the data set, with a consequent reduction of the computational effort. Moreover, the scheme is explicit.\\ 
Numerical tests show the accuracy and robustness of our approach to reconstruct curves and surfaces from relatively sparse data sets.
\end{abstract}

\section{Introduction}

Curvature-related level set models have been proposed in many problems of image and surface processing, due to their well understood analytical features and nice smoothing properties. In this paper, we consider one of such models, related to the problem of reconstructing a surface from a discrete set of points on it.

For $n\in\{2,3\}$, let $S=\{x_1,\ldots,x_N\}$ be a discrete set of points in $\R^n$, which are to be understood as points on a surface (or a curve if $n=2$), from which the surface itself has to be reconstructed. The set $S$ will be also termed as ``data set'' in what follows, and we assume that its points might be affected by noise. Thus, the reconstructed surface is not expected to pass through each of the data set points, but rather (depending on the criterion used to perform the reconstruction) to provide some trade-off between an exact interpolation and a smooth behaviour.

A level set model for this problem has been proposed by Zhao et al. in \cite{ZOMK}, and leads to the following evolutive problem
\begin{equation}\label{RicEq}
\begin{cases}
u_t (x,t)=d (x)\>{\rm{div}} \left(\frac{D u(x,t)}{|D u(x,t)|}\right)|D u(x,t)|+D d(x)\cdot Du(x,t)   \\
 u(x,0)=u^0(x) 
\end{cases}
\end{equation}
where $d(x)=d(x,S)$ is the euclidean distance from the set $S$, and $D$ denotes the gradient operator. As customary in level set methods, the reconstructed surface at a given time $t$ is represented as the zero-level set of the solution $u(x,t)$, that is,
\begin{equation}\label{gamma_t}
\Gamma_t = \{x\in\R^n: u(x,t)=0\}.
\end{equation}
In fact, \eqref{RicEq} is related to an energy functional in which the $L^1$ norm of the distance from $S$ is integrated on the whole surface (see \cite{ZOMK}). More precisely, given a surface $\Gamma$, we define the energy $\mathcal{E}$ as the surface integral
\[
\mathcal{E}(\Gamma) = \int_\Gamma d(x) ds,
\]
and look for a minimum of this functional (this has the clear meaning of a compromise between the total surface of $\Gamma$ and its distance from the data set $S$).
Once we define the evolution of an initial guess $\Gamma_0$ for the surface along the gradient flow of $\mathcal{E}$, and express the surface at time $t$ as the zero-level set of a function $u(x,t)$ as in \eqref{gamma_t}, we obtain \eqref{RicEq}. Details about the derivation of this model may be found in \cite{ZOMK}, whereas well-posedness of \eqref{RicEq} can be proved in the framework of viscosity solutions (see \cite{CIL92}), which requires minimal regularity assumptions on the solution.

Here, we are interested in the stationary version of \eqref{RicEq}, that is,
\begin{equation}\label{RicEqStaz}
d (x)\>{\rm{div}} \left(\frac{D u(x,t)}{|D u(x,t)|}\right)|D u(x,t)|+D d(x)\cdot Du(x,t) = 0,
\end{equation}
which plays the role of an Euler--Lagrange equation for the energy $\mathcal{E}$ and should be satisfied at local minima of the energy. In particular, a solution of \eqref{RicEqStaz} will be obtained in what follows as a regime solution of \eqref{RicEq} for $t\to\infty$.

We quote that similar techniques, retaining the regularizing effect of curvature-like terms, but possibly based on different evolution operators, are proposed for example in \cite{FMS, Li_etal}.

Of course, level set methods are not the only strategy used to solve the problem of surface reconstruction -- for example, segmentation techniques have been successfully proposed for this problem in recent years (see \cite{GBO09, LPZ13}). On the other hand, the largest amount of literature on the topic is probably devoted to least squares techniques. While we give up a complete review of this line of research, we mention that a relatively recent and successful technique in this framework makes use of Radial Basis Functions (RBF) space reconstructions, whose application to surface reconstruction stems from pioneering works published in the late 90s \cite{SPOK95, CFB97, Cetal01} (see also \cite{L02, W05} for a general review).

The aim of the present work is to investigate the application of semi-Lagrangian (SL) numerical techniques to \eqref{RicEq}, focusing in particular on their implementation with RBF space reconstructions. The application of RBF techniques to SL schemes has gained a certain popularity, although it has been restricted so far to the case of hyperbolic problems (see, e.g., \cite{BI02} and the literature therein).
On the other hand, Semi-Lagrangian schemes for curvature-related equations have been first proposed in \cite{FF03} and has gone through a number of improvements and applications (see, in particular, \cite{CFF10} for an in-depth convergence analysis and \cite{CF12, CF13} for two applications to image processing, along with \cite{FF13} for a general review on SL schemes). In this framework, RBF would be expected to provide a more flexible tool to construct a sparse space reconstruction. In fact, since their construction is not based on a space grid, Radial Basis Functions allow in principle for a sparse implementation, as well as for local refinements, although we are not aware of general strategies which could effectively handle RBFs in very disordered geometries. We propose therefore a {\it structured, but local} RBF space reconstruction in which one could better focus on the region close to the data set, instead of working on a full (and computationally more expensive) grid. This will be the final goal of the paper. We mention that, despite being constructed by different numerical tools, our scheme implements a ``localization'' of the numerical effort, much in the same spirit of the multigrid/multilevel techniques shown in \cite{BI02, Li_etal, Ye_etal}.

We finally remark that, while a model like \eqref{RicEq} might not be rated as the cutting-edge technique for surface reconstruction, yet the coupling of SL schemes with localized RBF space reconstructions shows a good potential in terms of accuracy and computational cost, and seems to be valuably applicable to a wider class of problems (in particular, level set models as the one under consideration).

\medskip

The outline of the paper is as follows. Section \ref{schema} will review the basic principles of construction of a SL scheme for \eqref{RicEq}, as well as the underlying ideas of RBF interpolation. Section \ref{test} will show numerical tests of increasing difficulty in two and three space dimensions, with both full and reduced grids.

\section{Numerical scheme}\label{schema}

We sketch in this section the construction of a SL approximation to \eqref{RicEq}. We start by sketching the basic ideas on the two dimensional case, then give the three-dimensional version of the scheme, and last describe the main improvements and modifications for the case of RBF space interpolations.

The main feature of the SL scheme under consideration is to be explicit, yet not constrained by the classic ``parabolic CFL'' condition, typical in the explicit treatment of diffusion terms. The degenerate diffusion performed by the curvature operator is treated by means of a convex combination of (interpolated) values of the numerical solution at the previous step, as we will soon show.

\subsection{The Mean Curvature operator}
Before introducing the scheme, we rewrite the mean curvature operator in such a way that the derivation of the method will be more natural.
 Let us recall that
\begin{eqnarray*}
 |Du(x,t)|{\rm{div}}\left(\frac{Du(x,t)}{|Du(x,t)|}\right) & = & \mathrm{tr}\left(P(Du)D^2u\right) = \\
 & = & {\rm{tr}}\left(\left(I-\frac{Du\otimes Du }{|Du|^2}\right)D^2u\right)
\end{eqnarray*}
where $I$ is the identity matrix, and
$$
P(Du) = I-\frac{Du\otimes Du }{|Du|^2}
$$
is a matrix which projects the diffusion on the tangent plane of each level surface. The projection $P(Du)$ is a matrix of rank $n-1$ with $n-1$ eigenvectors corresponding to the eigenvalue $\lambda=1$, and can be written as
\[
P(Du) = \sigma(Du)\sigma(Du)^\top,
\]
for a $n\times (n-1)$ matrix $\sigma(Du)$ having these eigenvectors as columns.

In the case $n=2$, there is only one eigenvector orthogonal to the gradient, namely
$$
\sigma(Du)= \frac{1}{|Du|}  \begin{pmatrix}
u_{x_2}\\ -u_{x_1} \\
\end{pmatrix}.
$$
In this case, the operator can be rewritten as
\begin{equation}\label{MCMop2d}
|Du(x,t)|\>{\rm{div}}\left(\frac{Du(x,t)}{|Du(x,t)|}\right)=\sigma(Du)^\top D^2 u\>\sigma(Du).
\end{equation}
In the case $n=3$, the projection matrix  $P$ is a $3\times 3$ matrix of rank 2 spanning the two-dimensional space orthogonal to the gradient of $u$. The two  orthonormal eigenvectors of $P$ are:
\begin{equation}\label{autovettori}
  \nu_1(Du)=\begin{pmatrix}
  \frac{-u_{x_3}}{\sqrt{u^2_{x_1}+u^2_{x_3}}} \\ \\
  0 \\ \\
  \frac{ u_{x_1}}{\sqrt{u^2_{x_1}+u^2_{x_3}}}
  \end{pmatrix}, \quad
  \nu_2(Du)=\frac{1}{|Du|}\begin{pmatrix}
  \frac{-u_{x_1}u_{x_2}}{\sqrt{u^2_{x_1}+u^2_{x_3}}} \\ \\
                        \sqrt{u^2_{x_1}+u^2_{x_3}} \\ \\
  \frac{-u_{x_2}u_{x_3}}{\sqrt{u^2_{x_1}+u^2_{x_3}}}
  \end{pmatrix},
\end{equation}
and, once set $\sigma(Du) =(\nu_1(Du),\nu_2(Du))$, the mean curvature  operator can be rewritten as
\begin{eqnarray}\label{MCMop3d}
|Du(x,t)|\>{\rm{div}}\left(\frac{Du(x,t)}{|Du(x,t)|}\right) & = &\mathrm{tr}\left(\sigma(Du)\sigma(Du)^\top D^2 u\right) = \nonumber \\
& = &  \frac{1}{2}(\nu_1+\nu_2)^\top D^2u\>(\nu_1+\nu_2)+\frac{1}{2}(\nu_1-\nu_2)^\top D^2u\>(\nu_1-\nu_2).
\end{eqnarray}
We will derive the scheme from this form, which corresponds to the probabilistic interpretation described in \cite{CFF10}.

\subsection{The basic 2D form}

In the two-dimensional case, the scheme has the structure
\begin{equation}\label{Scheme}
\begin{cases}
\displaystyle u^{n+1}_j = \frac{1}{2}I[u^n]\left(x_j+\D t\> D d(x_j)+\sqrt{2 \D t\> d(x_j) }\>\sigma^n_j\right)+ \\
\displaystyle \hspace{1.3cm} + \frac{1}{2}I[u^n]\left(x_j+\D t\> D d(x_j)-\sqrt{2 \D t\> d(x_j)}\>\sigma^n_j\right)  &x_j \in \G, n\ge 0 \\
 u^0_j = u^0(x_j) &x_j \in \G
 \end{cases}
\end{equation}
where $I[u^n](x)$ denotes a numerical interpolation of the solution $u^n$ computed at the point $x$, $\G$ denotes the set of grid nodes, and $\sigma^n_j$ is an approximation of the unit normal to $Du(t_n)$ defined by
\begin{equation}\label{sigma}
\sigma _j^n= \frac{1}{|D_j^n|}
\begin{pmatrix}
D_{2,j}^n \\ -D_{1,j }^n \\
\end{pmatrix},
\end{equation}
for some consistent approximations $D_{1,j }^n$, $D_{2,j }^n$ (e.g., centered partial incremental ratios of the numerical solution) for the partial derivatives $u_{x_1}(x_j,t_n)$, $u_{x_2}(x_j,t_n)$, with $|D_j^n|$ denoting the euclidean norm of the approximate gradient at $(x_j,t_n)$, i.e.,
$$
|D_j^n| = \left({D_{1,j }^n}^2 + {D_{2,j }^n}^2\right)^{1/2}.
$$
Note that, at each node, \eqref{Scheme} performs an average of the solution at the two points $x_j+\D t\> D d(x_j)\pm\sqrt{2 \D t\> d(x_j)}\sigma^n_j$. The term $\D t\> D d(x_j)$ clearly represents an upwinding along the advection field $D d(x_j)$, as usual in SL schemes. On the other hand, the degenerate operator related to curvature is treated by adding a further displacement $\pm\sqrt{2 \D t\> d(x_j)}\>\sigma^n_j$, which has the effect of generating a diffusion along the tangent space of the level sets, according to the curvature operator \eqref{MCMop2d}. More explicitly, \eqref{Scheme} can be rewritten via minor algebraic manipulations as
\begin{eqnarray}\label{forma2}
\frac{u^{n+1}_j-u^n_j}{\D t} & = &  \frac{d(x_j)}{|h_j|^2}\Big(I[u^n](x_j+\D t\> D d(x_j)+h_j) - 2 I[u^n](x_j+\D t\> D d(x_j)) + \nonumber \\
&& + I[u^n](x_j+\D t\> D d(x_j)-h_j)\Big) + \frac{1}{\D t} \Big(I[u^n](x_j+\D t\> D d(x_j)) - u^n_j\Big)
\end{eqnarray}
where $h_j=\sqrt{2 \D t\> d(x_j)}\>\sigma^n_j$ is a finite increment along the tangent space (approximately) spanned by $\sigma^n_j$. In the right-hand side of \eqref{forma2}, it is easy to recognize that the two terms are consistent with respectively the second and the first directional derivatives which appear in \eqref{RicEq}. A more detailed consistency analysis for the pure curvature case is carried out in \cite{FF03, CFF10}.

Singular points at which $D_j^n\approx 0$ (where the definition \eqref{sigma} would not make sense) are to be treated by suitable techniques. A typical choice (see \cite{CFF10}) is to replace the degenerate diffusion with an isotropic one (e.g., the heat operator) whenever the approximate gradient is below a given threshold, say, $|D_j^n|< C\D t^\alpha$, for a proper choice of the constant $C$ and the exponent $\alpha>0$. More explicitly, in a form parallel to \eqref{forma2}, by computing the diffusion term as a 5-point laplacian we obtain
\begin{eqnarray}\label{forma2_sing}
\frac{u^{n+1}_j-u^n_j}{\D t} & = &  \frac{d(x_j)}{2|h_j|^2}\Big(I[u^n](x_j+\D t\> D d(x_j)+|h_j|e_1) + I[u^n](x_j+\D t\> D d(x_j)-|h_j|e_1) + \nonumber \\
&& + I[u^n](x_j+\D t\> D d(x_j)+|h_j|e_2) + I[u^n](x_j+\D t\> D d(x_j)-|h_j|e_2) - \nonumber \\
&&  -4 I[u^n](x_j+\D t\> D d(x_j)) \Big) + \frac{1}{\D t} \Big(I[u^n](x_j+\D t\> D d(x_j)) - u^n_j\Big)
\end{eqnarray}
(where $e_1$ and $e_2$ are the base vectors of $\R^2$), which results in a convex combination of the values $I[u^n](x_j+\D t\> D d(x_j) \pm |h_j|e_i)$ and is therefore explicit and stable.
This technique of treating singularities can be shown to be compatible with the definition of viscosity solution for \eqref{RicEq} (see \cite{CFF10}). 

\subsection{The 3D case}

Following \cite{CFF10}, we can also write a three-dimensional version of the scheme (which will be used in the section on numerical tests to recover surfaces). The new form of $\sigma_j^n$ will then be written as
\begin{equation*}
\sigma_j^n =
\begin{cases}
\left(\nu_1(D_j^n),\nu_2(D_j^n)\right) & \text{ if }\sqrt{{D_{1,j }^n}^2+{D_{3,j }^n}^2}\neq 0 \\
(d_1,d_2) & \text{ if } \sqrt{{D_{1,j }^n}^2+{D_{3,j }^n}^2} = 0,
\end{cases}
\end{equation*}
where $d_1=(1,0,0)^\top$ and $d_2=(0,0,1)^\top$. Using this matrix to perform an average of 4 points on the tangent plane, we write the scheme as
\begin{equation}\label{Scheme3d}
u^{n+1}_j = \frac{1}{4} \sum_{i=1}^4 I[u^n]\left(x_j+\D t\> D d(x_j)+\sqrt{2 \D t\> d(x_j) }\>\sigma^n_j\delta_i\right)
\end{equation}
in which the vectors $\delta_i$ $(i=1,\ldots,4)$ are defined as
\[
\delta_i = \begin{pmatrix} \pm 1 \\ \pm 1 \end{pmatrix}
\]
for all combinations of the signs.
Via minor algebraic manipulations, \eqref{Scheme3d} can be rewritten as
\begin{eqnarray}\label{forma3}
&&\hspace{-.8cm}\frac{u^{n+1}_j-u^n_j}{\D t}  = \frac{1}{\D t} \Big(I[u^n](x_j+\D t\> D d(x_j)) - u^n_j\Big)+\\
&& + \frac{d(x_j)}{|h^1_j|^2}\Big(I[u^n](x_j+\D t\> D d(x_j)+h^1_j) - 2 I[u^n](x_j+\D t\> D d(x_j)) + I[u^n](x_j+\D t\> D d(x_j)-h^1_j)\Big)  +\nonumber \\
 && + \frac{d(x_j)}{|h^2_j|^2} \Big(I[u^n](x_j+\D t\> D d(x_j)+h^2_j) - 2 I[u^n](x_j+\D t\> D d(x_j)) + I[u^n](x_j+\D t\> D d(x_j)-h^2_j))\Big). \nonumber 
\end{eqnarray}
where
$$
h^1_j=\sqrt{2 \D t\> d(x_j)}\>\left(\nu_1(D_j^n)+\nu_2(D_j^n)\right)
$$
$$
h^2_j=\sqrt{2 \D t\> d(x_j)}\>\left(\nu_1(D_j^n)-\nu_2(D_j^n)\right).
$$
In the second and third line of \eqref{forma3}, it is possible to recognize the second finite differences along the directions $\nu_1(D u)+\nu_2(D u)$ and $\nu_1(D u)-\nu_2(D u)$,  which has the effect of generating a diffusion along the tangent space of the level sets, in agreement with the curvature operator \eqref{MCMop3d}. In the first  line  of \eqref{forma3}, an upwind approximation of the transport term appears.

In (approximately) singular conditions, i.e., when $|D_j^n|<C\D t^\alpha$, the diffusion term is switched to a 7-point laplacian by analogy with \eqref{forma2_sing}.

\subsection{RBF implementation}

While the SL scheme \eqref{Scheme} has proved to be robust and relatively accurate in a variety of applications, we study in this work an adaptation to this specific case. In \eqref{RicEq}, the interest is in following the zero-level set of the solution, which is in turn supposed to stay in the neighbourhood of the data set $S$. The computational effort has therefore to be concentrated in this latter region -- computing an accurate solution away from the data set is useless.

Keeping this idea in mind, we implement \eqref{Scheme} with a space reconstruction in the form of a Radial Basis Function (RBF) interpolation, which lends itself to a sparse implementation (see \cite{L02}). In particular, we have used here the Matlab RBF interpolation toolbox described in \cite{C06}.

The general structure of the space reconstruction under consideration is
\begin{equation}\label{rbf}
I[u](x) = c_0(u) + c(u)\cdot x + \sum_i \lambda_i(u) \phi_\rho(|x-x_i|)
\end{equation}
where the scalar $c_0$, the vector $c$ and the coefficients $\lambda_i$ (all of which depend on $u$) are determined by imposing interpolation conditions, and a suitable closure of the system (see \cite{C06}). Note that in this case, since $I[u]$ has a globally defined form, the computation of the $\sigma_j^n$ in \eqref{Scheme}--\eqref{Scheme3d} might be carried out explicitly.

In \eqref{rbf}, the function $\phi_\rho : \R^+\to \R$ provides the radial term. Although (as for the linear RBF defined by $\phi(t)=t$) this term might not depend on further parameters, it is usually expected to depend on a shape constant $\rho$ (e.g., the ``variance'' in a Gaussian function). Whenever $\rho$ appears in the definition of $\phi$, it is generally recognized that the choice of this parameter has a critical effect on the accuracy. Moreover, optimal strategies for selecting $\rho$ exist for a regular distribution of the nodes $x_i$, but not if the set $\G$ of nodes is genuinely scattered. On the other hand, a structured, full-grid implementation of RBF space reconstruction gives essentially no advantage with respect to a more conventional interpolation. For all these reasons, we will implement the following strategy:

\begin{itemize}

\item
In order to have a clear guideline for choosing the scale parameter $\rho$, nodes are distributed over a {\em reduced, but structured set}, i.e., we start from a structured orthogonal grid and keep only the nodes belonging to some given neigbourhood of the data set. For example, we can use nodes $x_j$ for which
\begin{equation}\label{n_band}
d(x_j,S) < \delta_S,
\end{equation}
for some suitable threshold $\delta_S$. The nodes of the data set are also used as RBF nodes. As an example, Fig. \ref{Test1dgrid} shows data set points as circles in red, and RBF nodes as asterisks in blue. Note that taking $\delta_S=O(k\D x)$ results in creating a reduced grid containing about $k$ layers of nodes both inside and outside $S$.

\item
It is known that RBF are particularly unsuitable to extrapolate a function at points far from the region covered by the nodes.
In order to keep the nodes concentrated near the data set, still avoiding instabilities when far from this region, we introduce a set $\mathcal{A}$ of {\em anchor nodes} in which the numerical solution has a prescribed value, possibly both inside and outside the curve/surface to be reconstructed. Anchor nodes essentially work as boundary conditions, by locking the value of the numerical solution at ``far'' points where the computation of RBF interpolation cannot be (and in fact is not required to be) accurate.

Referring again to Fig. \ref{Test1dgrid}, anchor nodes appear as the outer crosses in black. In our numerical tests, the anchor nodes have been placed on a frame surrounding the data set, following the geometry of the edges of the computational domain. In any of the tests, no anchor node has been used inside the data set.

\end{itemize}

We will show in the section on numerical examples that this framework results in a considerable reduction of the computational grid, still preserving robustness of the algorithm.

\paragraph{Initial condition}

In order to speed up convergence, the initial condition should have a zero-level set $\Gamma_0=\partial\Omega_0$ {\em as close as possible to the surface (or curve) to be reconstructed}, although, of course, the surface itself is unknown. Moreover, a natural choice for the initial condition $u^0$ would be the signed distance from the initial interface $\Gamma_0$, that is,
\begin{equation*}
d_s(x,\Gamma_0) :=
\begin{cases}
d(x,\Omega_0) & \text{if } x\notin \Omega_0 \\
-d(x,\Omega_0) & \text{if } x\in \Omega_0.
\end{cases}
\end{equation*}
In practice, since we look for an asymptotic solution for \eqref{RicEq}, it turns out that neither the choice of $\Gamma_0$ nor the structure of $u_0$ is critical -- clearly, a better initial conditions results in a faster solution. In lack of a sharper information, a relatively conservative choice of the initial condition might be set as
\begin{equation}\label{Inicond}
u_j^0 = \begin{cases}
|x_j|^2 - R^2 & \text{ if } x_j\in\mathcal{G} \\
C & \text{ if } x_j\in\mathcal{A}
\end{cases}
\end{equation}
with a suitable choice of the constants $R,C>0$ in order to start from a zero-level set enclosing the data set, but not too far from it. This choice, which will be used in the numerical examples of the next section, avoids the cost of computing the signed distance, but maintains negative values for $u$ inside the surface $\Gamma_t$, and positive values outside. We did not experience any critical behaviour of the scheme with respect to this choice, nor the need for a re-initialization of the level set function. In fact, the effect of the advection term in \eqref{RicEq} is to push the solution towards the data set, and this causes the solution to maintain a steep transition in the neighbourhood of the interface. 

\section{Numerical tests}\label{test}

We present in this section some numerical examples of shape reconstructions with the algorithm described above. In the two-dimensional case, we will also compare the reduced grid with the full grid scheme. Note that, in two dimensions, it is known (see \cite{ZOMK}) that a polygonal line joining the points of the data set is a local minimizer of the energy, so this should be regarded as an ``exact solution'' for the 2D case.
In all the test, we compute the normalized $L^1$-norm of the update between two successive iterates :
\begin{equation}\label{error}
E^n_1=\frac{\sum_{j\in \mathcal{G}} |u^{n+1}_j-u^n_{j}|}{\sum_{j\in \mathcal{G}} |u^n_j|},
\end{equation}
this being an indication on the speed of convergence to the regime solution. Note that, when plotted on a semi-log scale (as it will be done in what follows), a linear decrease of $E_1^n$ suggests a contraction-like behaviour of the scheme.

In the following examples, we will keep the data set in a quite sparse form. We found it discriminating to work with a relatively low number of points, and this will be our strategy in order to check the robustness of the numerical algorithm.

\subsection{Heart-shaped 2D data set}

\paragraph{Full grid}

\begin{figure}
\begin{center}
 \epsfig{figure=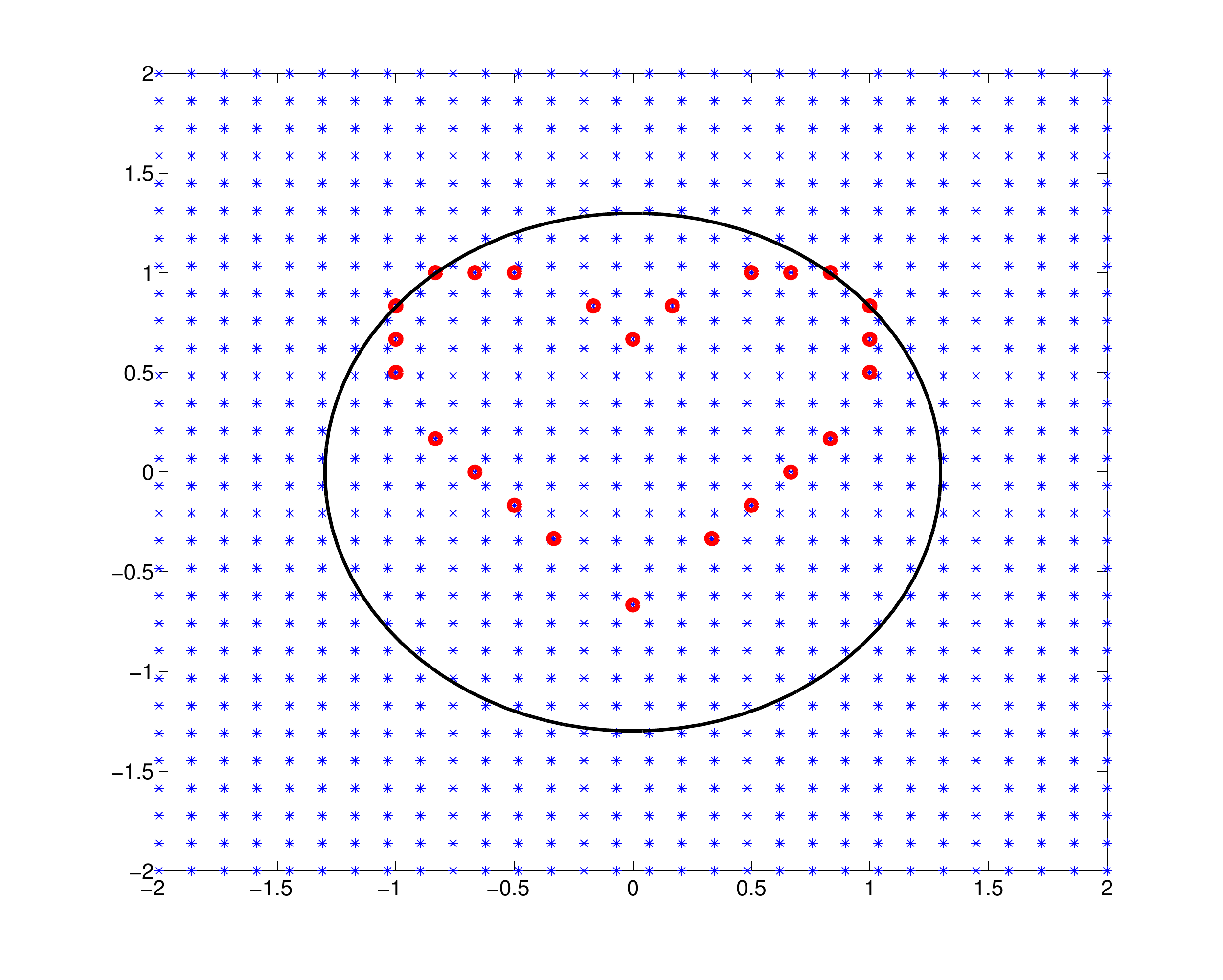,width=8cm}
 \epsfig{figure=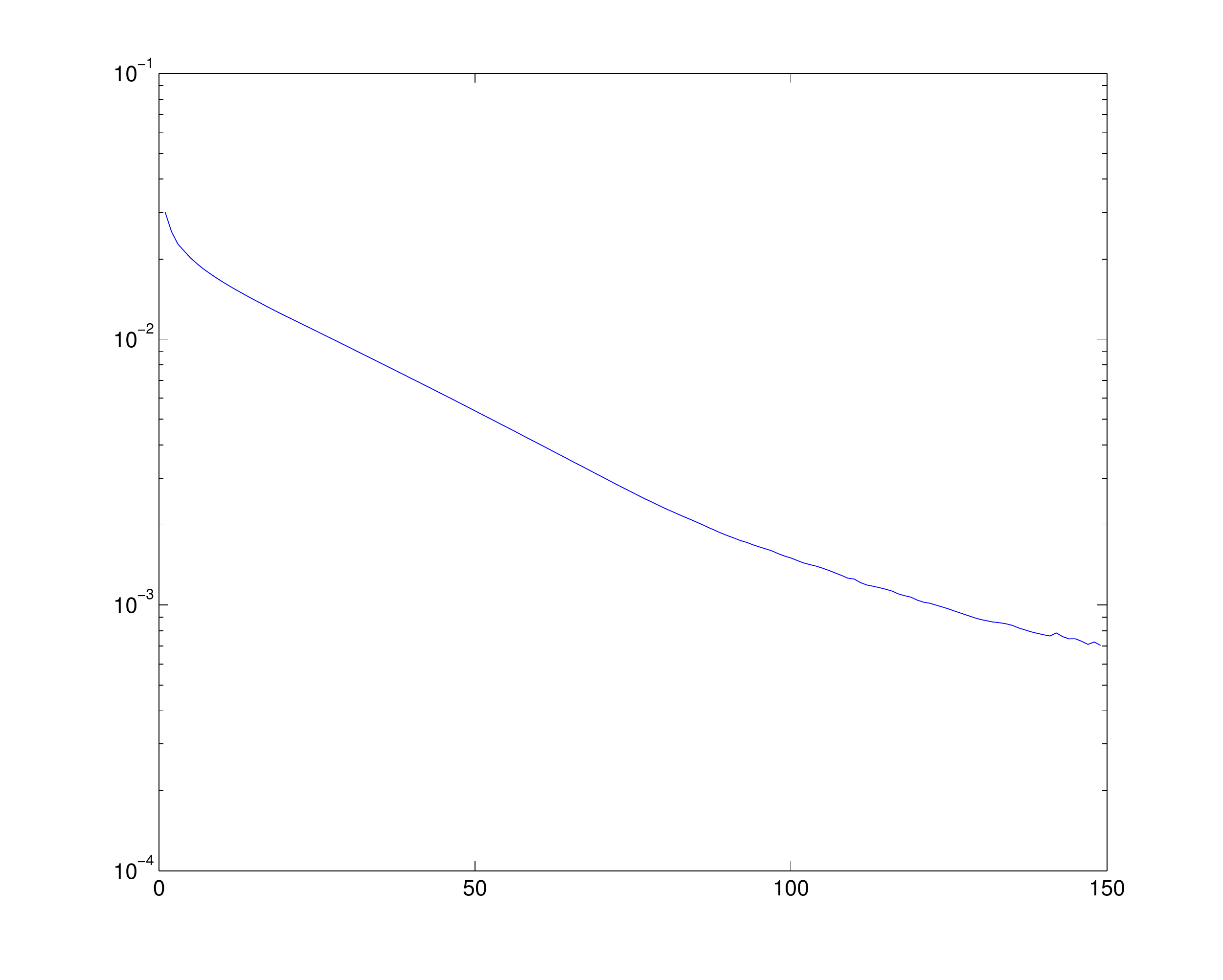,width=8cm}
\caption{{Left: Grid nodes $\G\setminus S$ (blue asterisks), data set $S$ (red circles) and initial condition (continuous line) . Right: relative update between two successive iterations.}}
\label{Test1bgrid}
\vspace{1cm}
\epsfig{figure=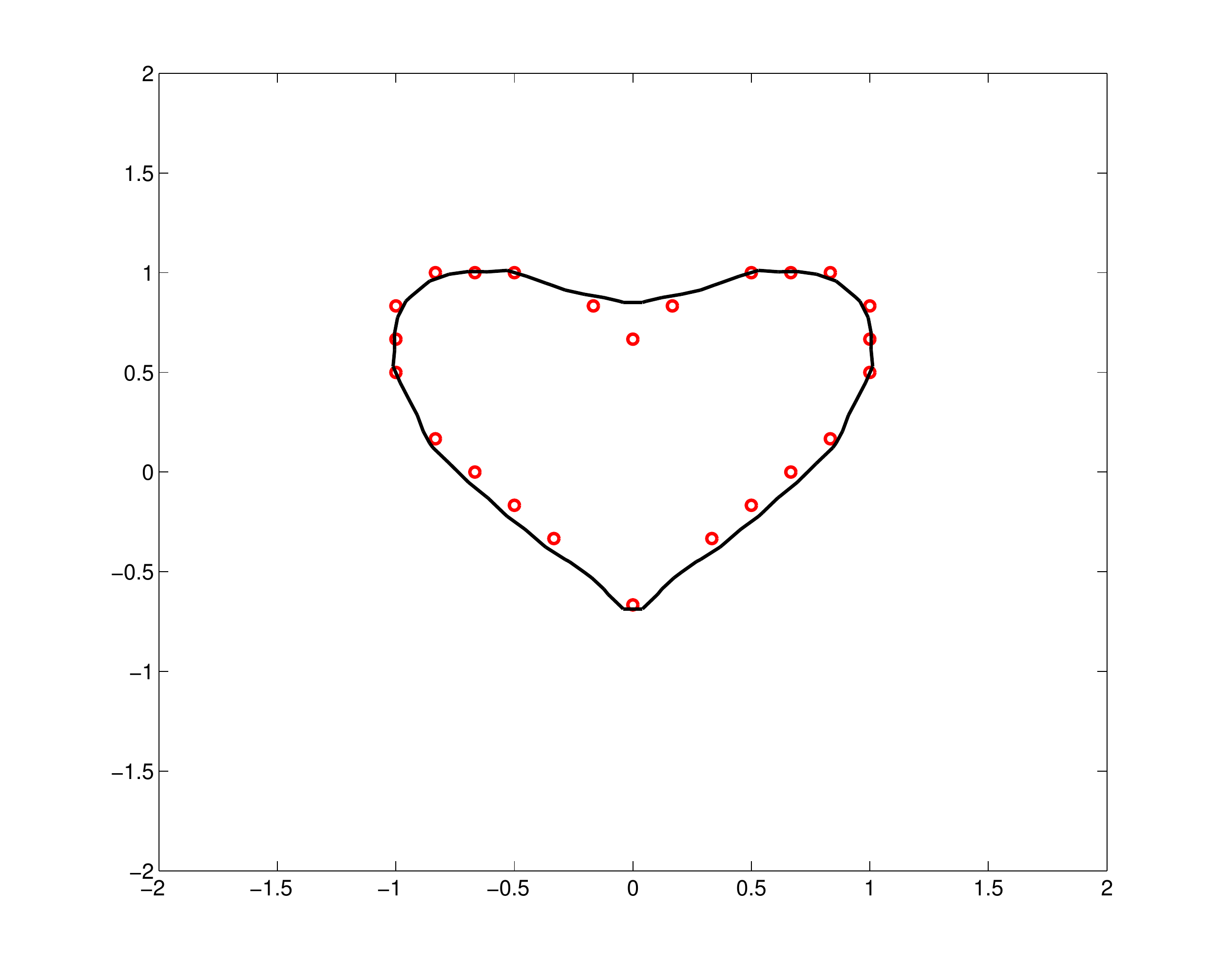,width=5cm}  
\epsfig{figure=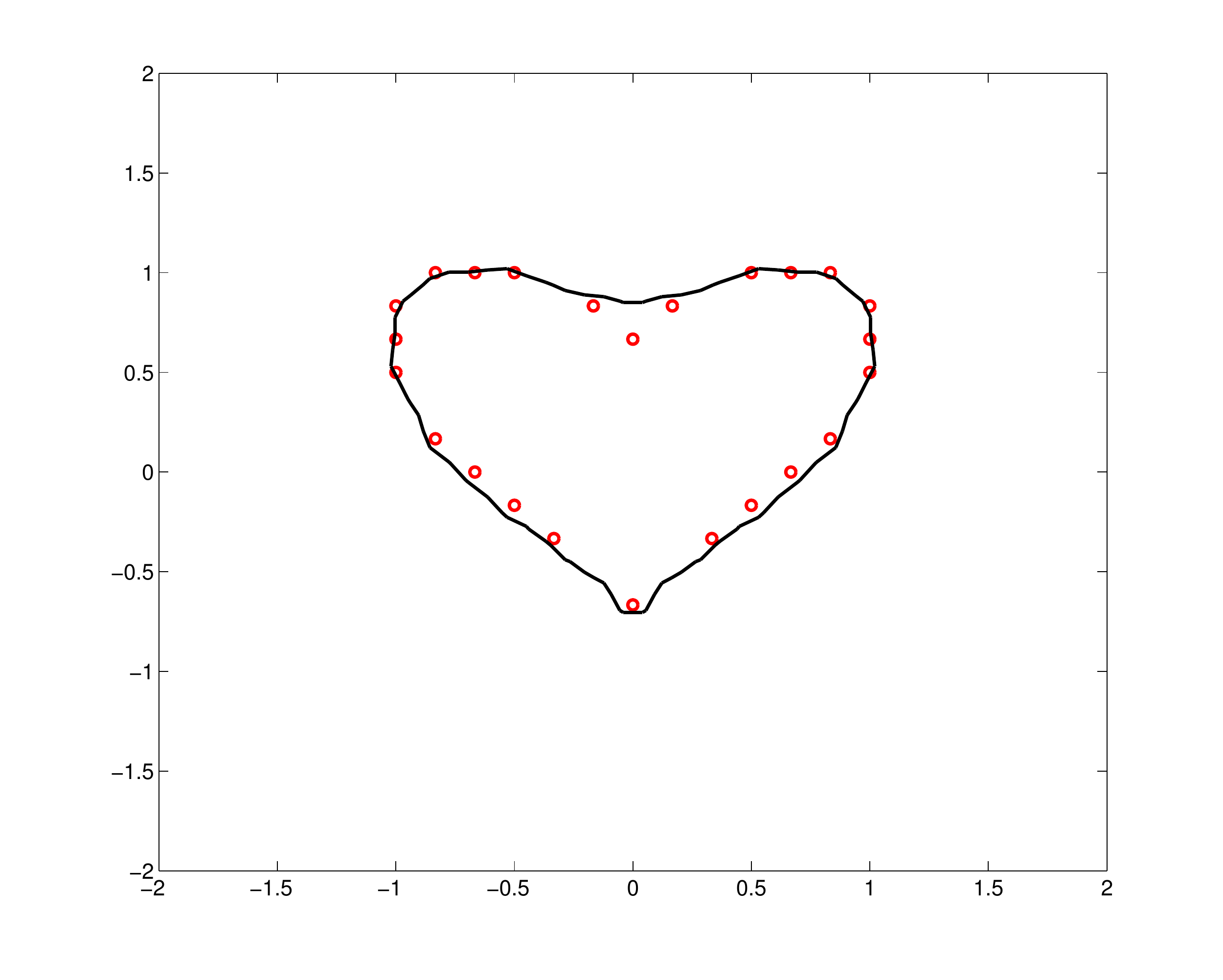,width=5cm}   
\epsfig{figure=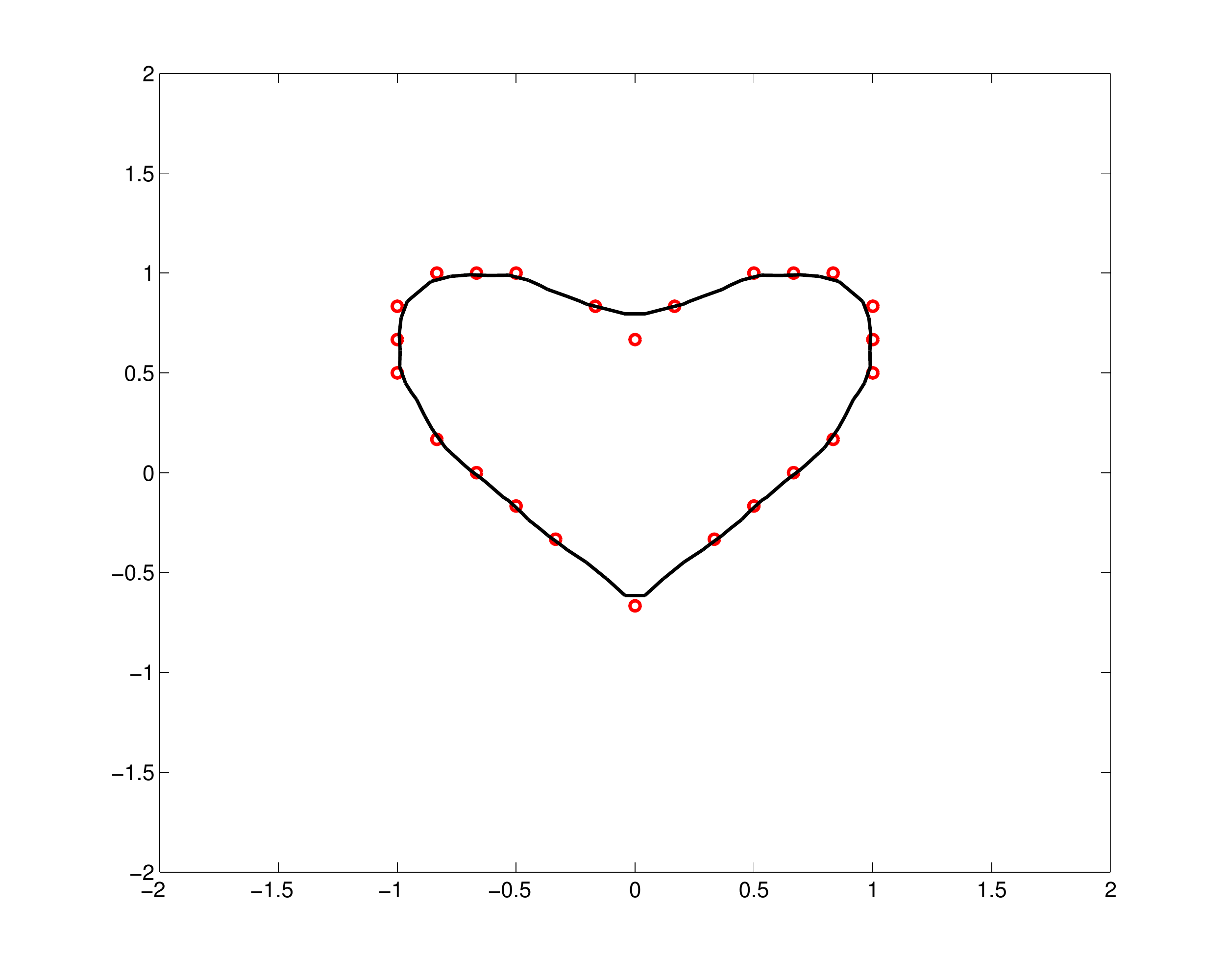,width=5cm}   
\caption{{Solution at the final time (continuous line) and data set (red circles).  Left: multiquadric RBFs. Center: linear RBFs. Right: linear RBFs with finer mesh.}}
\label{Test1bsol}
\end{center}
\end{figure}

We first consider a data set $S$ made of 24 points uniformly chosen on a heart-shaped curve and a grid $\G$ of $30\times30$ evenly spaced nodes on a domain $[-2,2]^2$, as shown in the left plot of Fig. \ref{Test1bgrid}. We apply \eqref{Scheme} on the grid $\G$ with time step $\D t=0.01$ and show (in the right plot of Fig. \ref{Test1bgrid}) the trend of the normalized update \eqref{error} between two successive iterations, the algorithm being stopped after 150 iterations. The final solution obtained by multiquadric RBFs (with a scale factor $\rho=\Delta x$, see \cite{C06}) is shown in the left plot  of Fig. \ref{Test1bsol}. The final solution obtained by linear RBFs and by linear RBFs applied on a finer $60\times60$ grid are shown respectively in the center and in the right plot  of Fig. \ref{Test1bsol}.
The reconstructed curve is very close to the theoretical forecast of a polygonal line, except for a slight smoothing of the non-convex, upper section of the shape (see Fig. \ref{Test1bsol}). Regardless of the number of points in the data set, this gives the natural indication that the computational grid is still too coarse. Note also that convergence of iterates tends to become quite slow, although it seems not to depend on the type of RBF used. In fact, we have included a single plot for both linear and multiquadric RBFs, whose convergence histories show an analogous behaviour.

\paragraph{Reduced grid}

\begin{figure}
\begin{center}
 \epsfig{figure=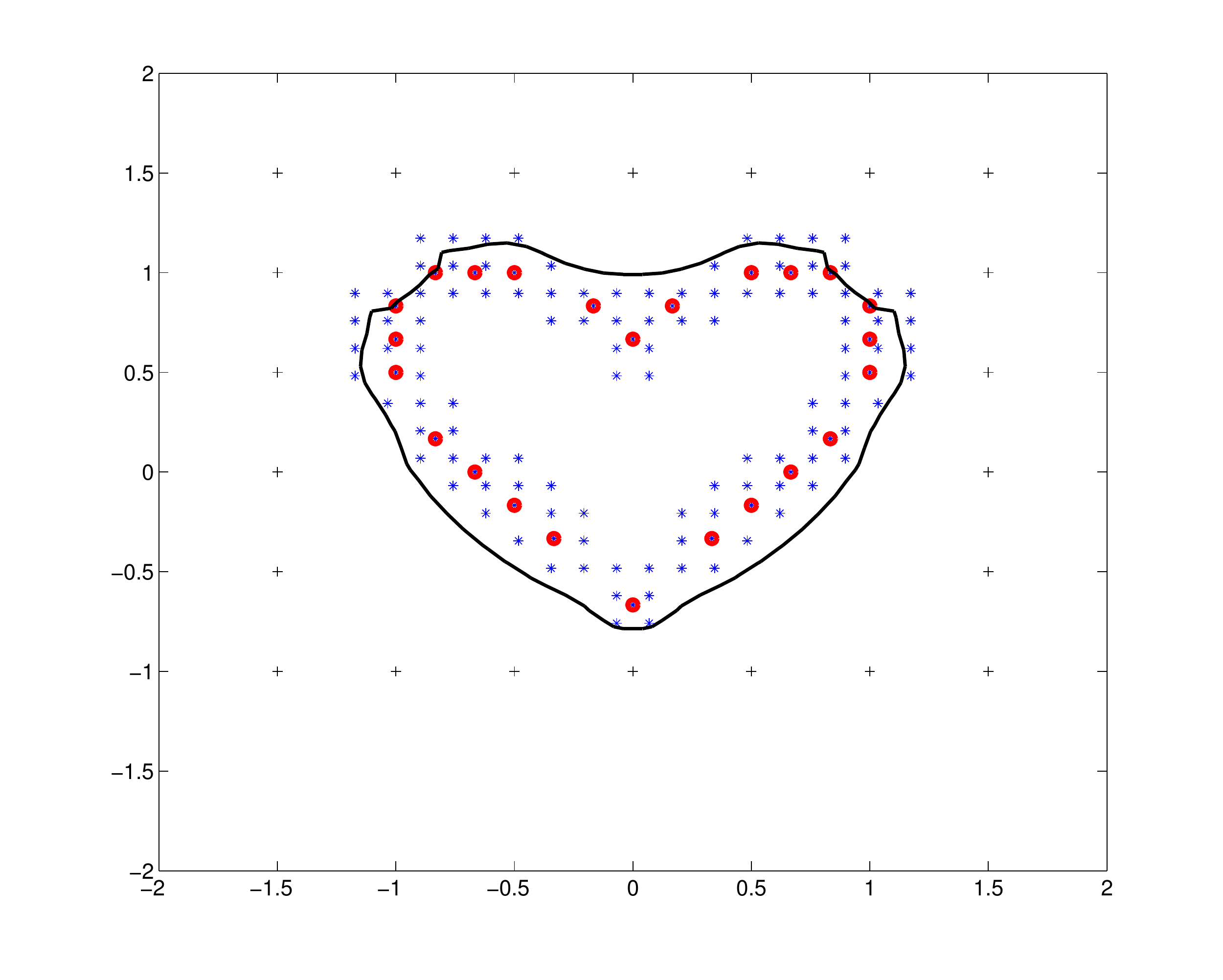,width=8cm}
 \epsfig{figure=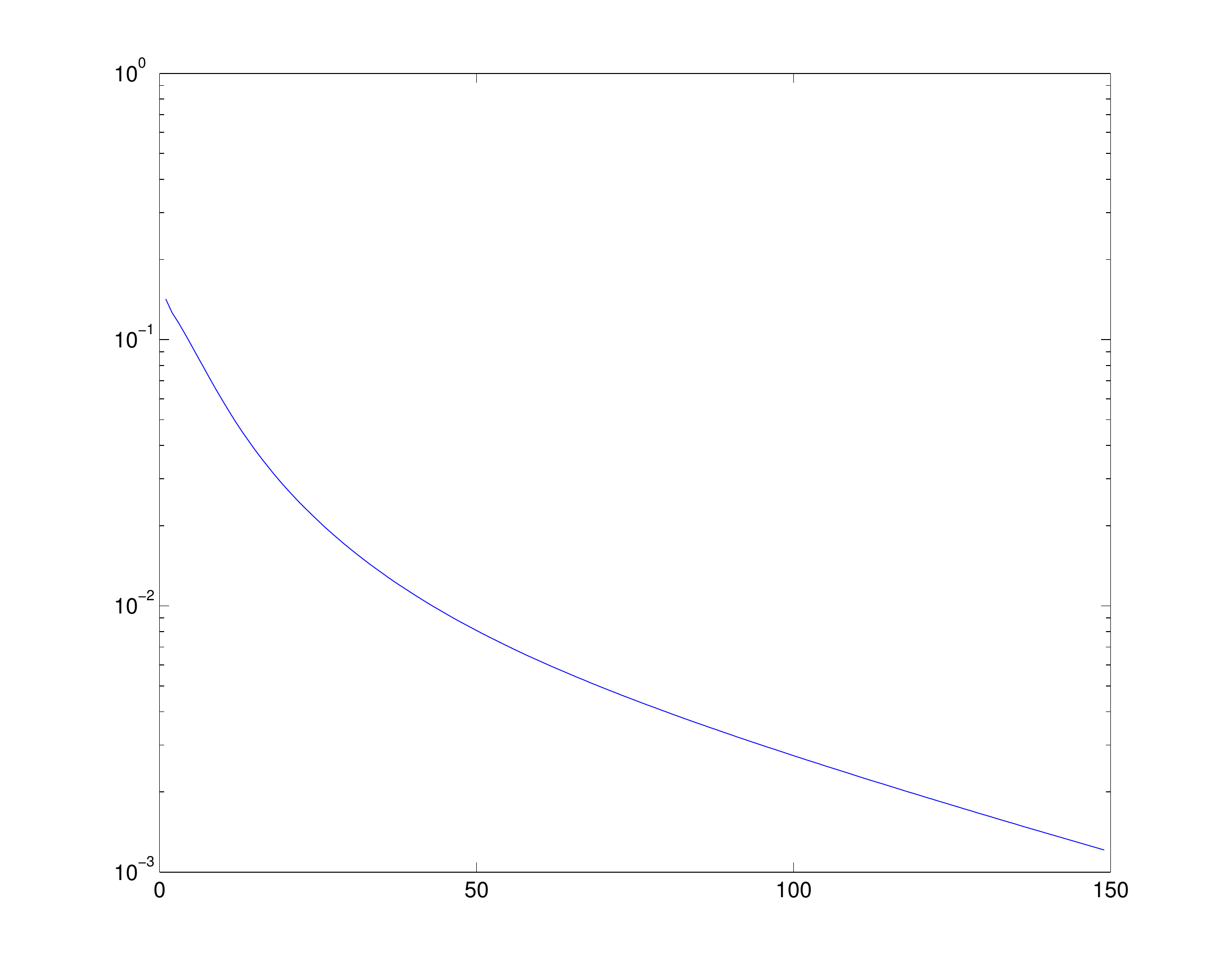,width=8cm}
\caption{{Left: Grid nodes $\G\setminus S$ (blue asterisks), data set $S$ (red circles) and initial condition (continuous line) . Right: relative update between two successive iterations.}}
\label{Test1dgrid}
\vspace{1cm}
\epsfig{figure=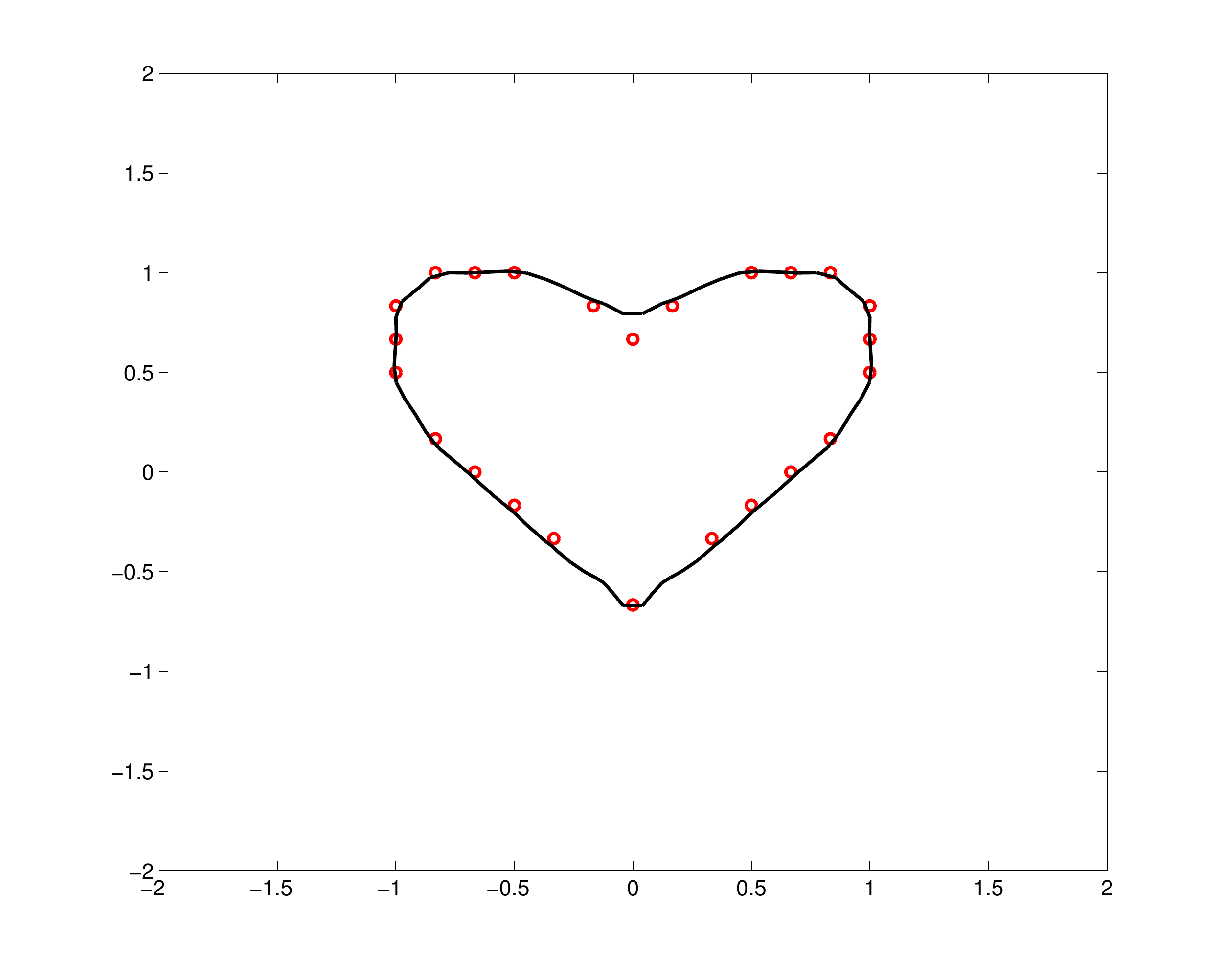,width=8cm}
\epsfig{figure=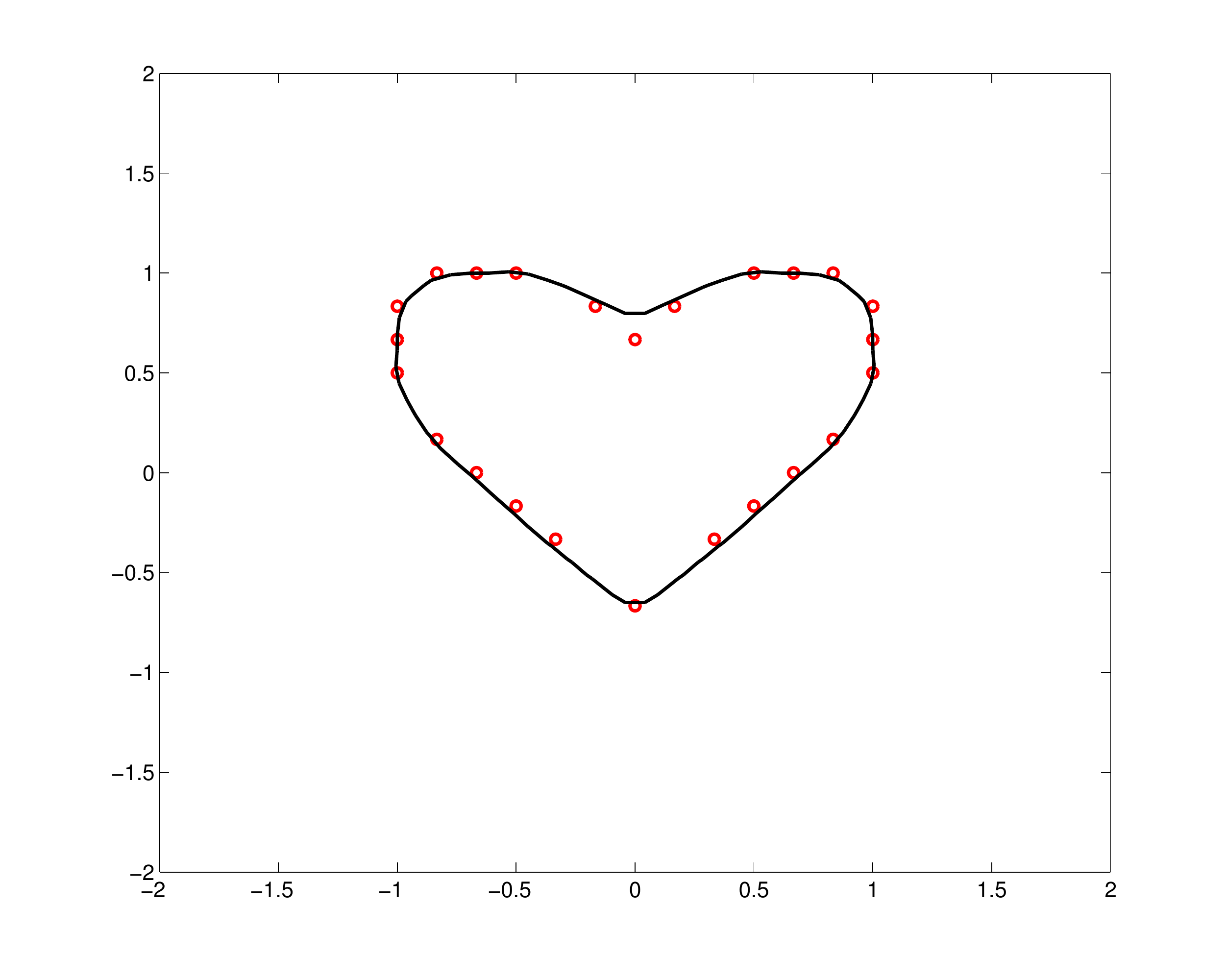,width=8cm}  
 \caption{{Solution at the final time (continuous line) and data set (red circles).  Left: multiquadric RBFs. Right: linear RBFs.}}
\label{Test1dsol}
\end{center}
\end{figure}

The same test is performed again on a reduced grid $\G$ made of 106 point distributed around the data set according to \eqref{n_band} with $\delta_S=0.2$, as shown in the left plot of Fig. \ref{Test1dgrid}, which also shows the anchor points (marked as the outer black crosses) together with the initial condition. The smaller area covered by the grid allows to decrease the total number of nodes to a fraction of about 12\% of the original grid. The final solution obtained by multiquadric RBFs with $\rho=\Delta x$ and by linear RBFs are respectively shown  in the left and right  plots of Fig. \ref{Test1dsol}.
The result is essentially equivalent to what is obtained with the full grid (Fig. \ref{Test1dsol}). Note that, comparing the right plots of Fig. \ref{Test1bgrid} and \ref{Test1dgrid}, convergence of the iterative solver seems also to be comparable in terms of iteration number, but with a lower cost for a single iteration. As in the case of a full grid, convergence of the iterates for linear and multiquadric RBFs is very similar.

\subsection{3D data sets}

\paragraph{Heart-shaped surface}

\begin{figure}
\centering
\epsfig{figure=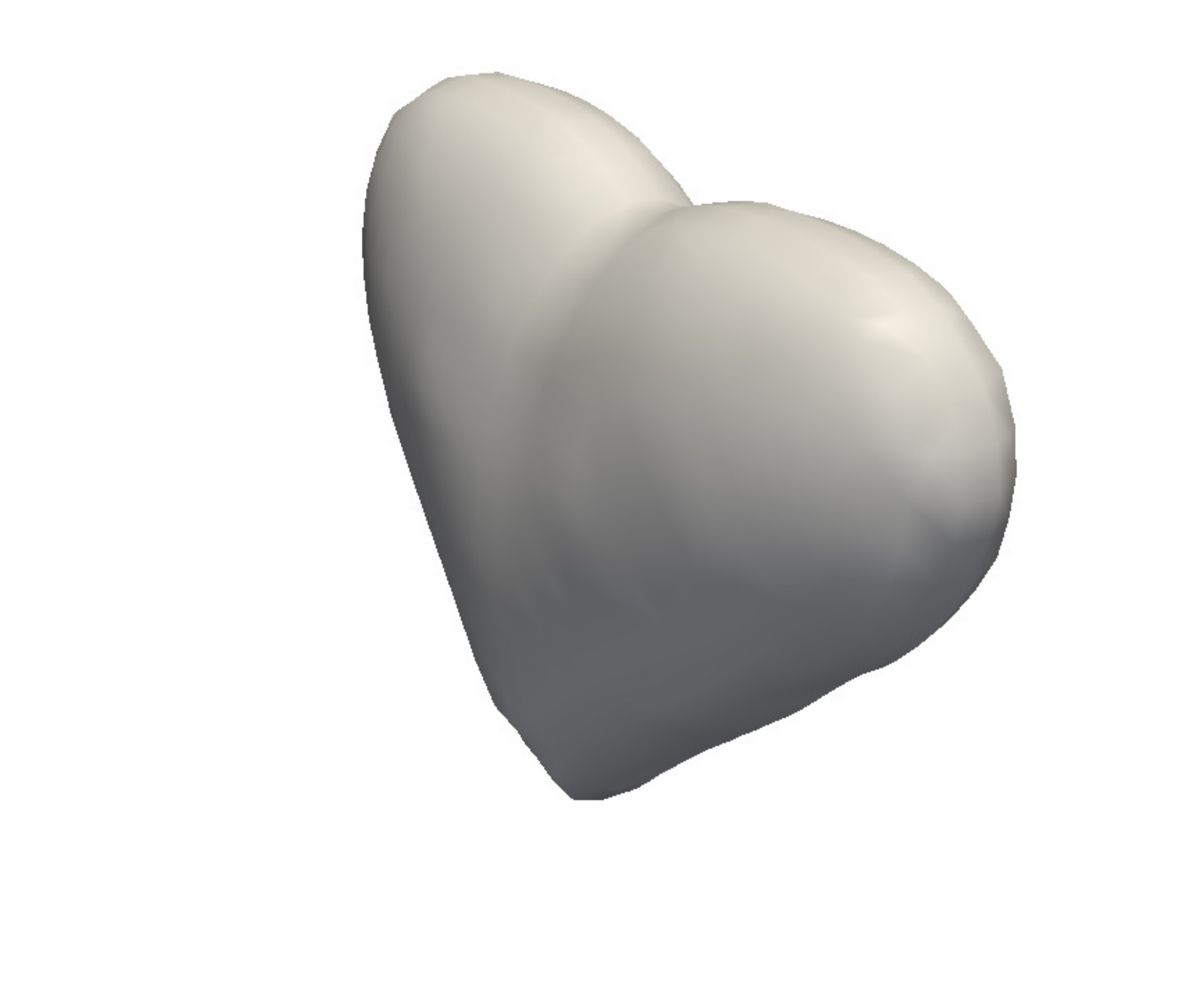,width=8cm}\epsfig{figure=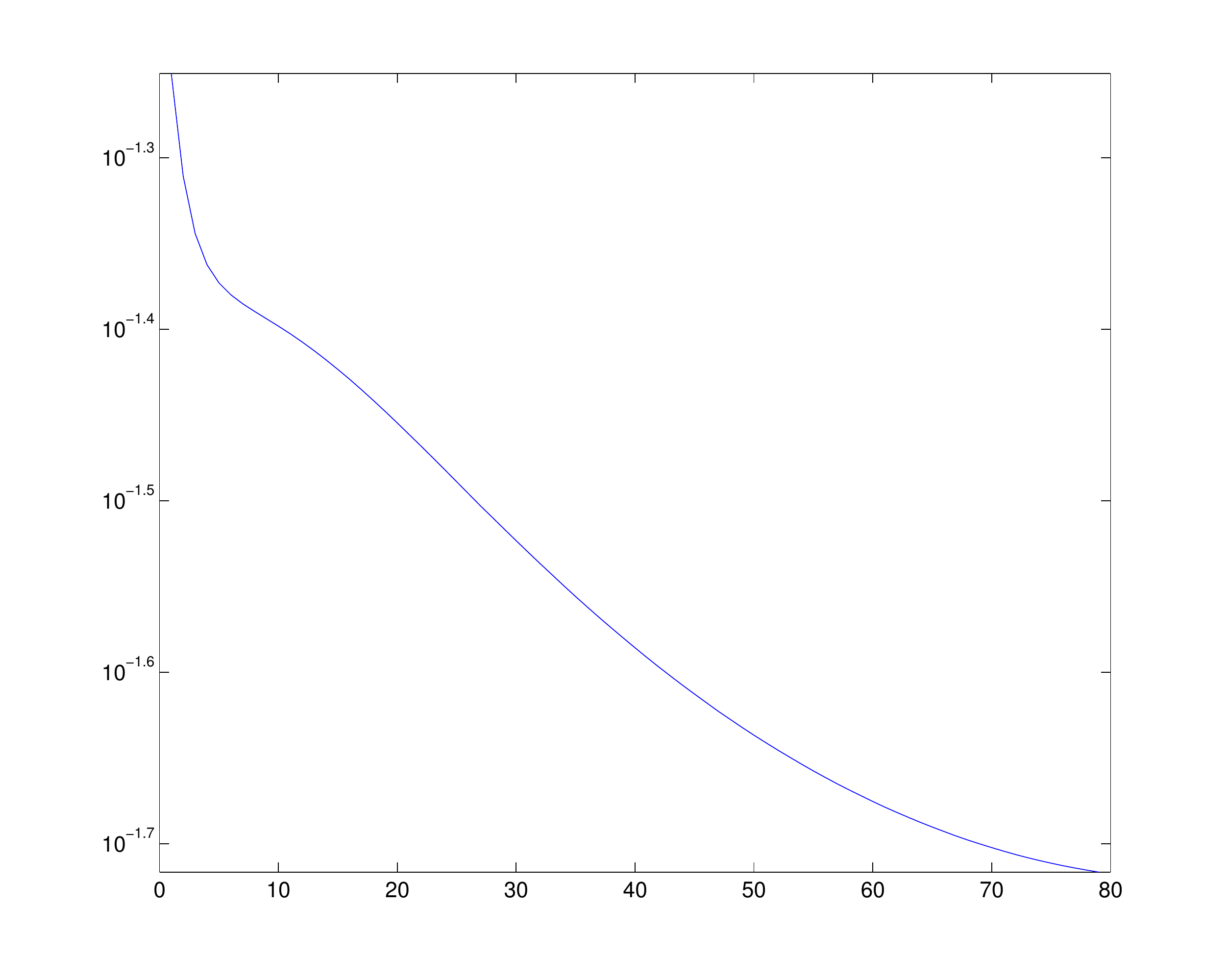,width=8cm}
\caption{{Left: reconstructed surface for the 3D heart shape. Right: relative update between two successive iterations}
\label{Test2bsol}}
\epsfig{figure=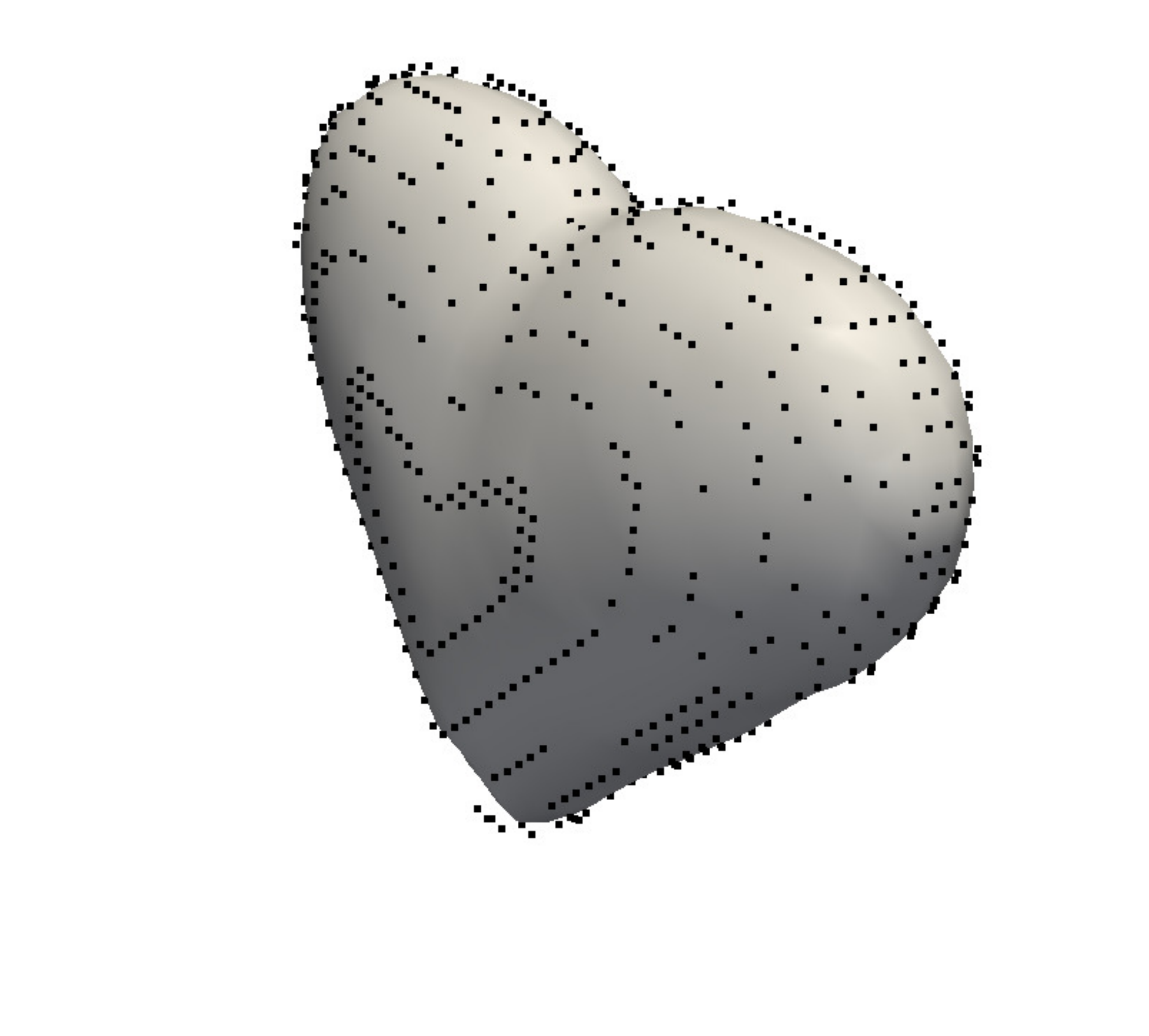,width=8cm}\epsfig{figure=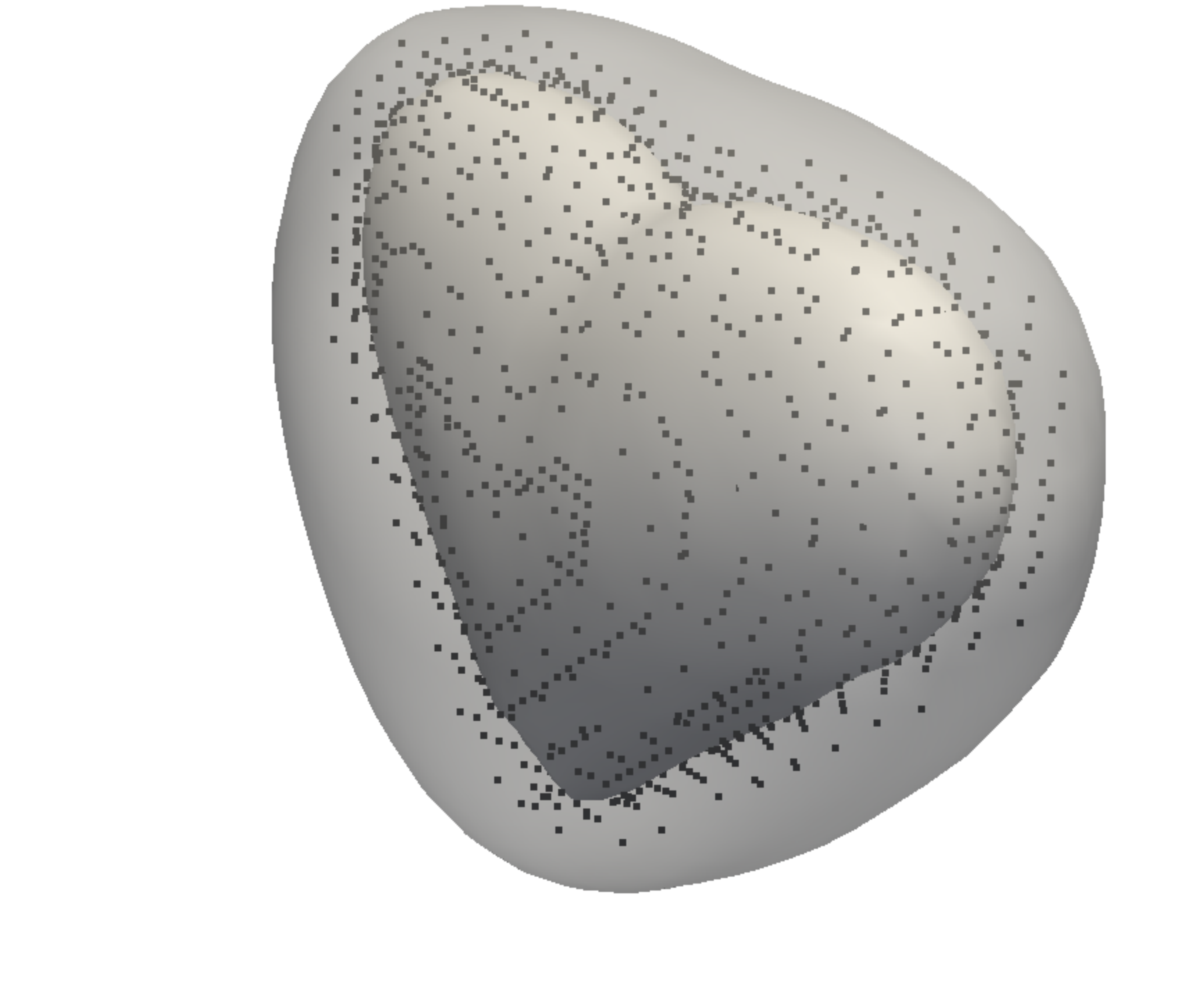,width=8cm}
\caption{{Left: reconstructed surface compared with the data set (black dots). Right: reconstructed surface compared with the grid (black dots) and the zero-level set of the initial condition.}\label{Test2bgrid}}
\epsfig{figure=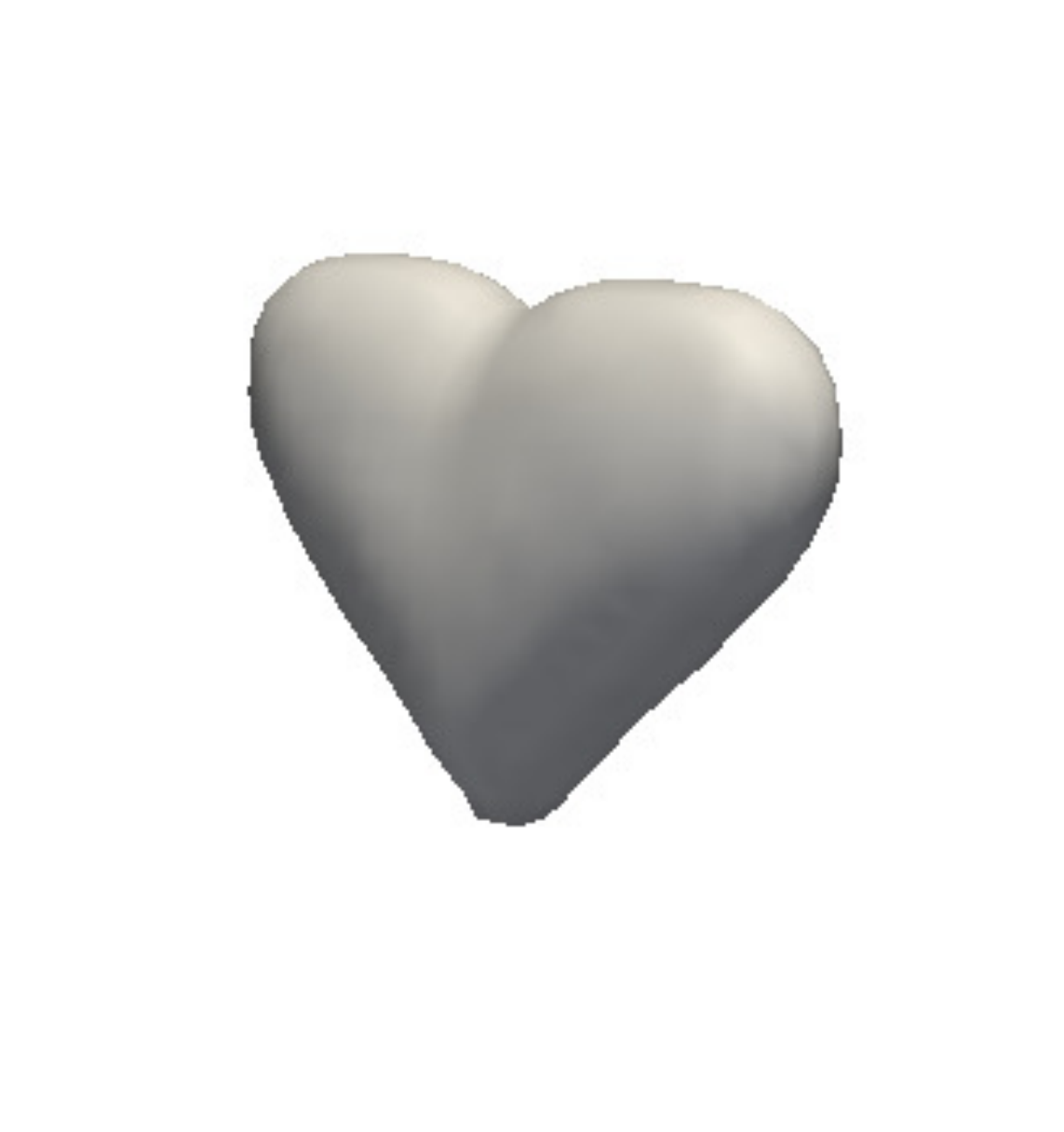,width=5.5cm}\epsfig{figure=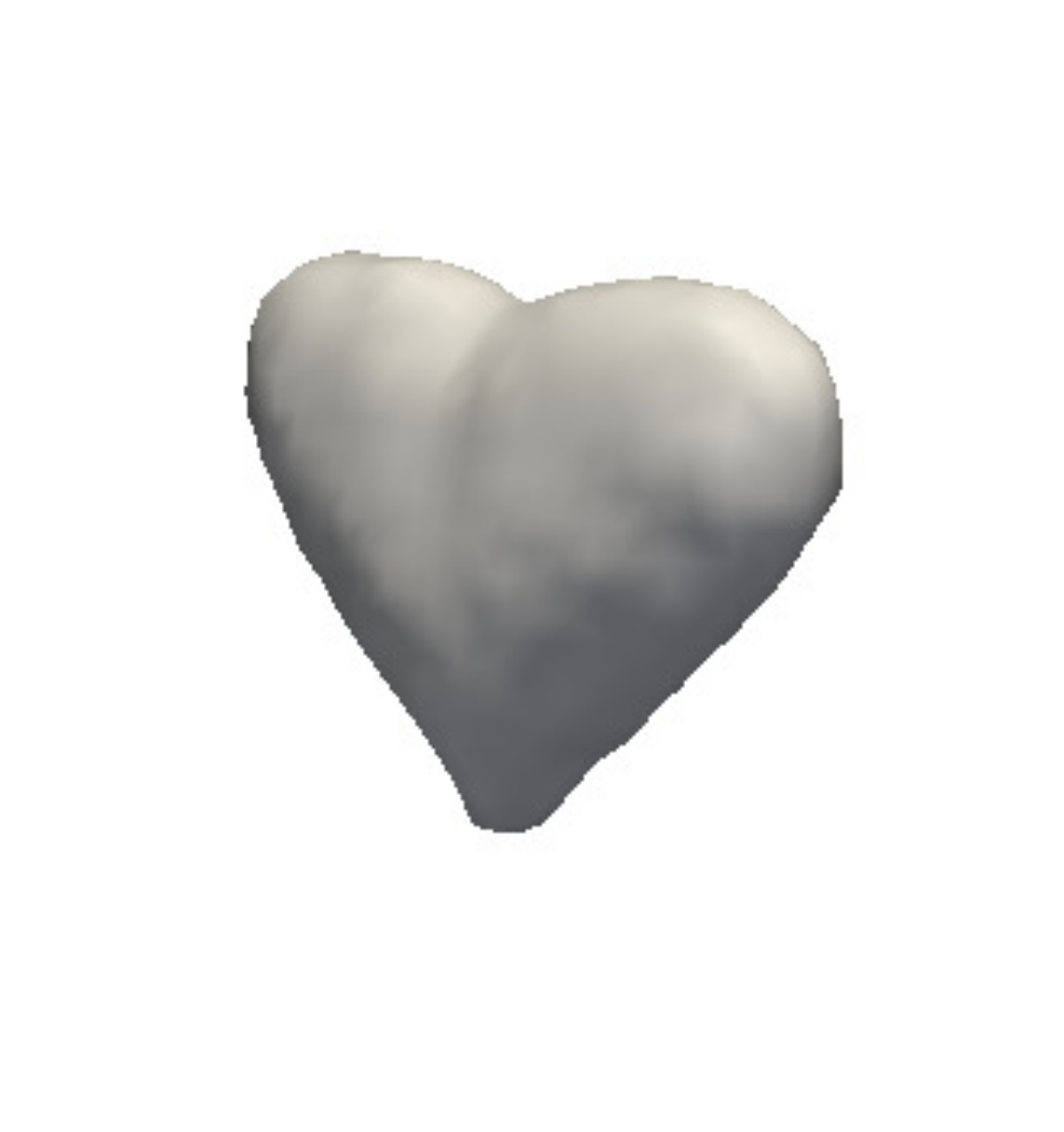,width=5.5cm}\epsfig{figure=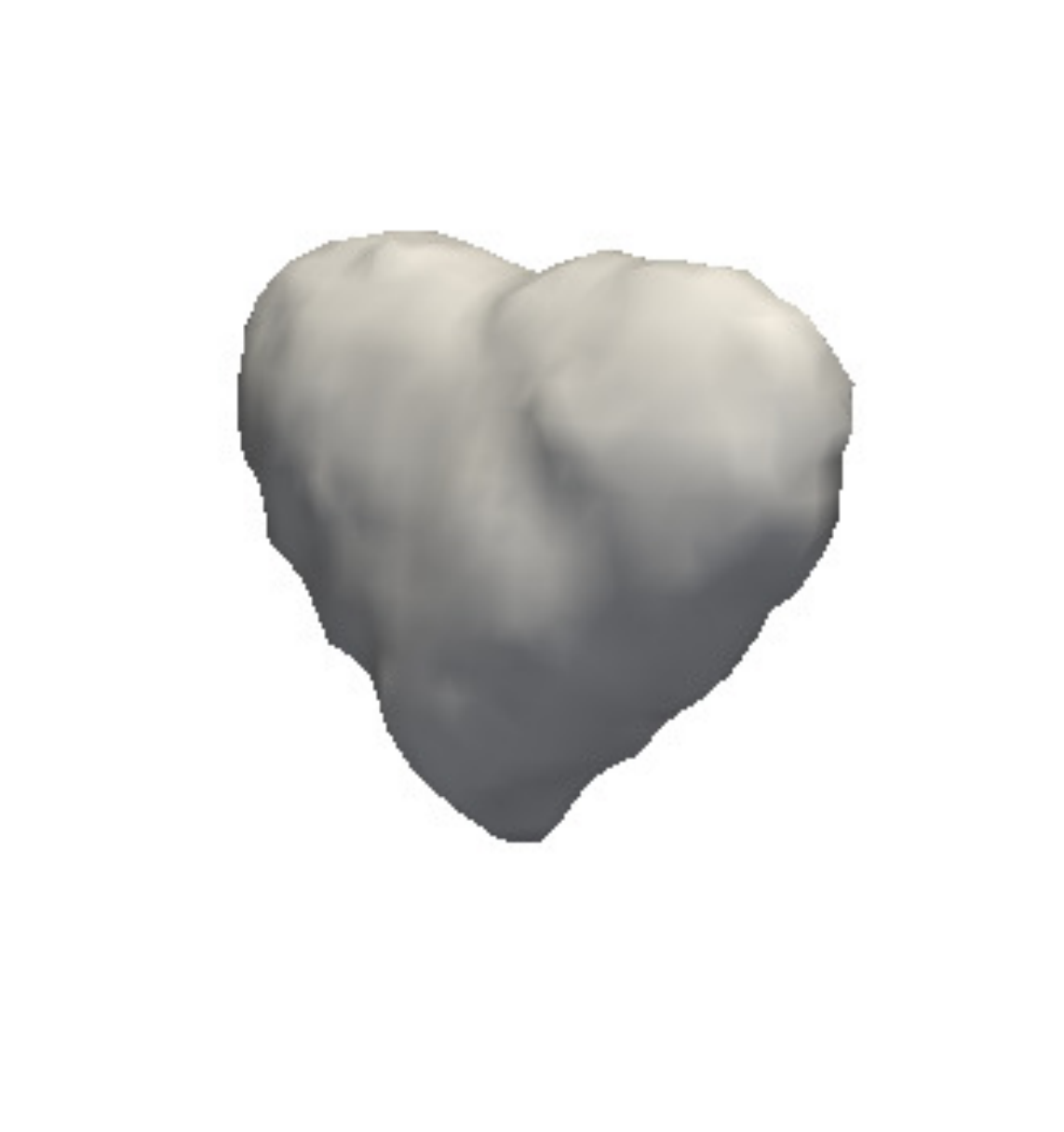,width=5.5cm}
\caption{{Reconstructed surface from a noisy data set, with $\eta=0.01$ (left), $\eta=0.025$ (center) and $\eta=0.05$ (right). }\label{Test2bnoise}}
\end{figure}

We  consider a data set $S$ made of 748 points chosen on a 3D heart shape in the set $[-2,2]^3$. We build a grid  $\G$ made by the dataset plus $945$ additional nodes uniformly distributed on a narrow band around the data set, computed according to \eqref{n_band} with $\delta_S=0.1$ from  a full grid of $40\times40\times40$ nodes.\\
We apply \eqref{Scheme3d} on the grid $\G$ with time step $\D t=0.005$. The right plot of Fig. \ref{Test2bsol} shows the norm of the update between two successive iterations. The algorithm is stopped after 80 iteration, and the final solution is shown in the left plot of Fig. \ref{Test2bsol}. Dataset, grid and zero-level set of the initial condition are shown in Fig. \ref{Test2bgrid}.

We also analyze in this example the robustness of the method with respect to noisy data. To this end, we perturb the data set points by adding random displacements, uniformly distributed in $[-\eta,\eta]^3$. The recovered surfaces obtained by applying \eqref{Scheme3d} to three different levels of noisy datasets ($\eta=0.01,0.025,0.05$), still on the same space grid $\G$, with time step $\D t=0.005$ and 80 iterations, are shown in Fig. \ref{Test2bnoise}.

\paragraph{Two intersecting cubes }

\begin{figure}
\centering
\epsfig{figure=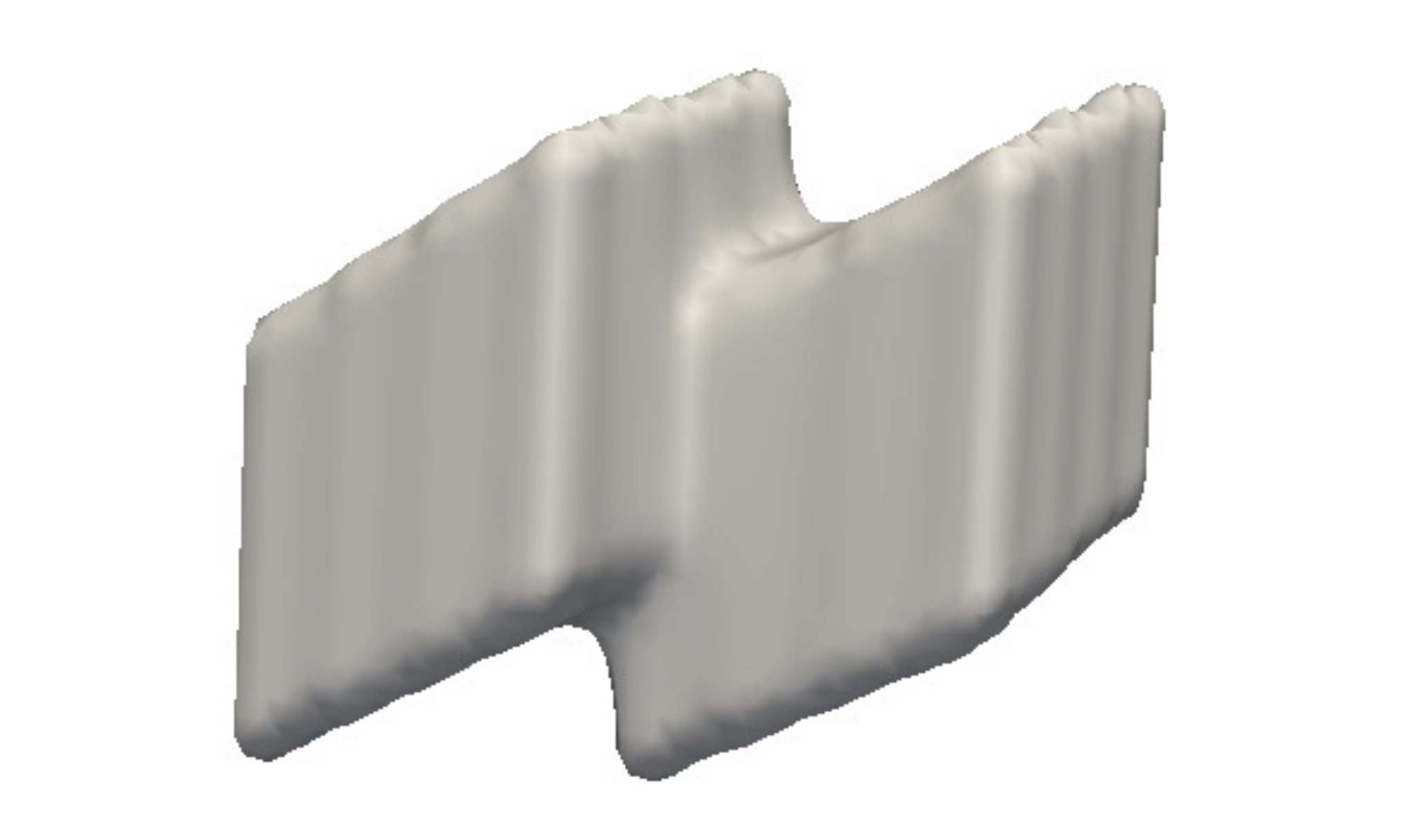,width=8cm}\epsfig{figure=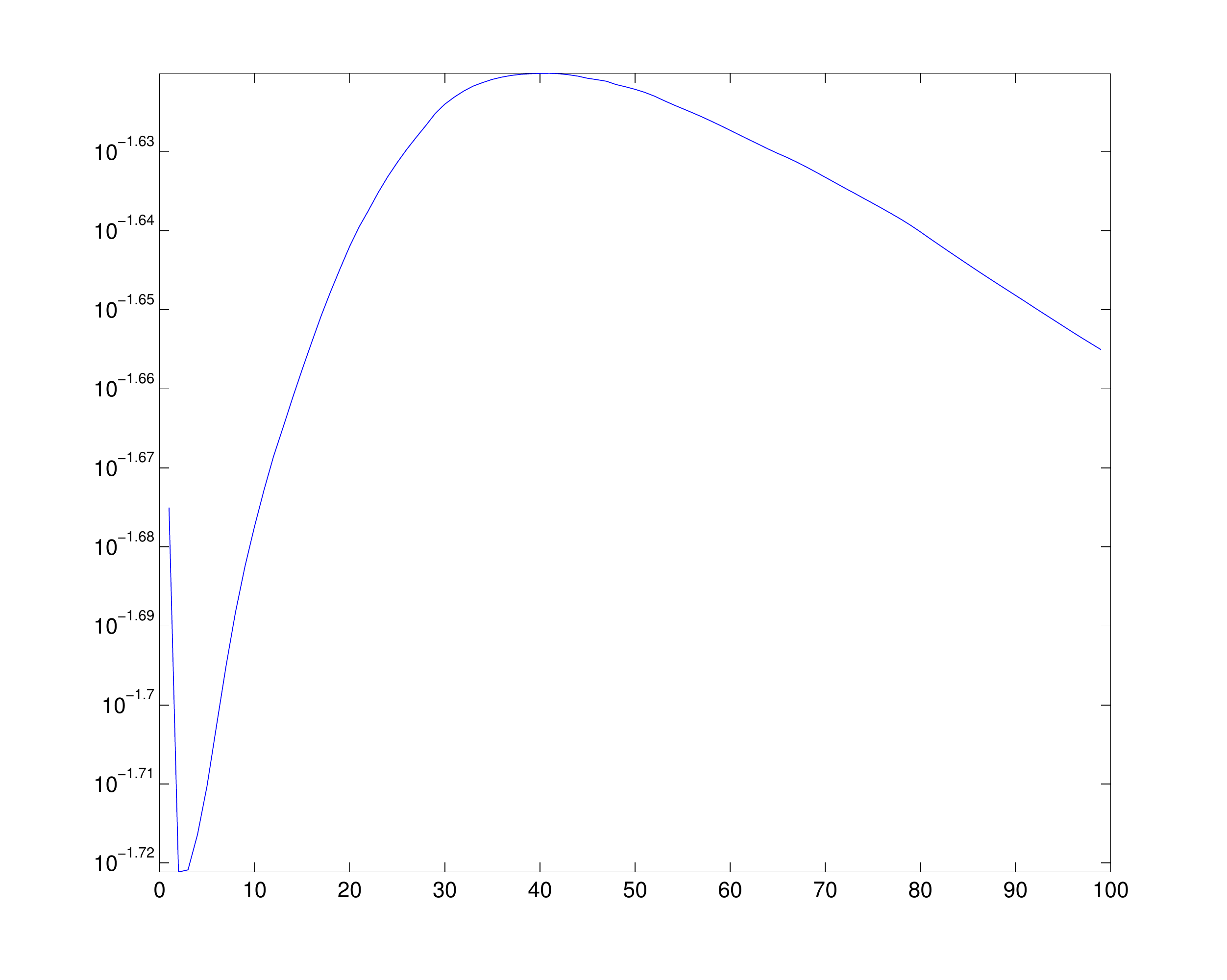,width=8cm}
\caption{{Left: reconstructed surface for the intersected cubes. Right: relative update between two successive iterations}
\label{cubisol}}
\epsfig{figure=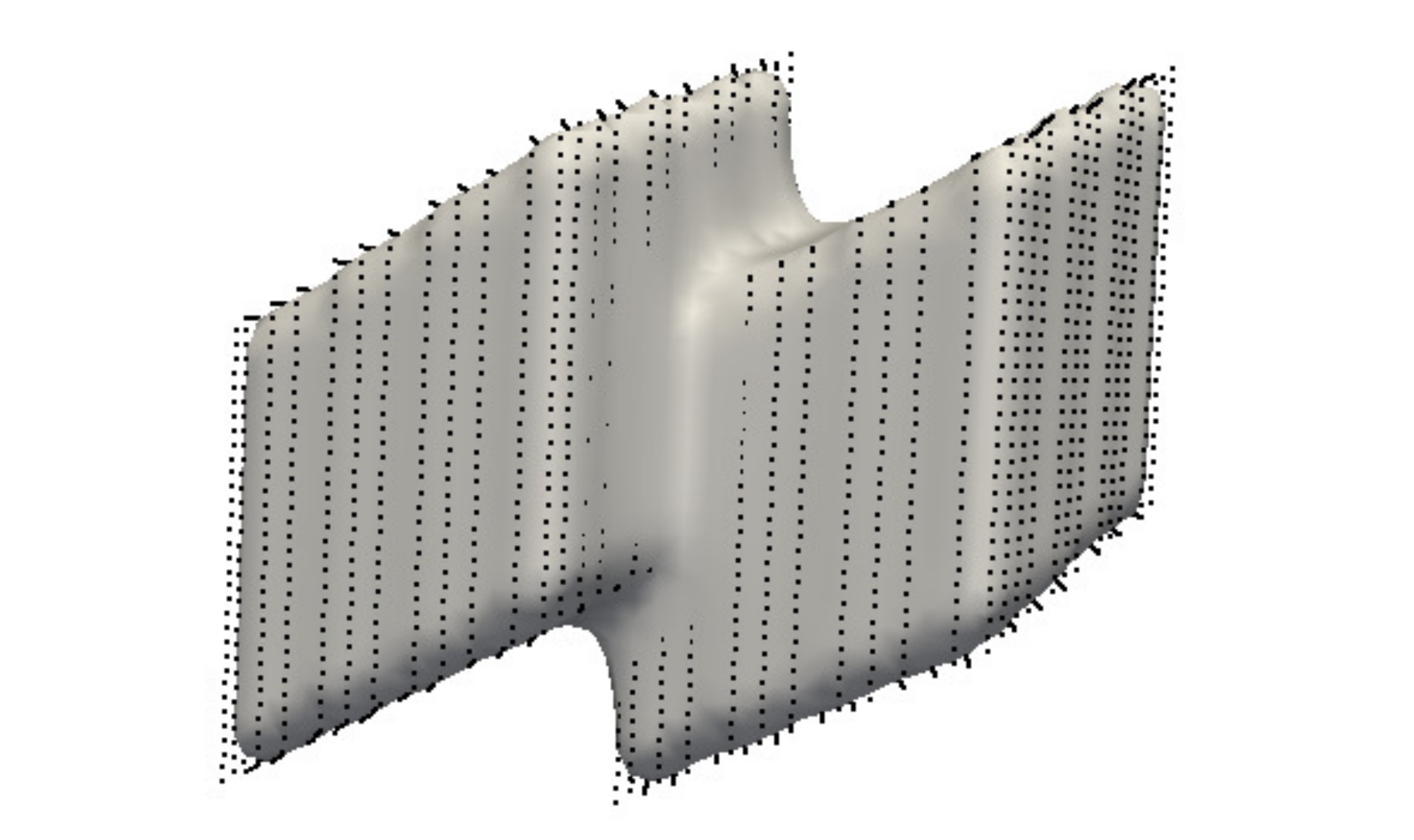,width=8cm}\epsfig{figure=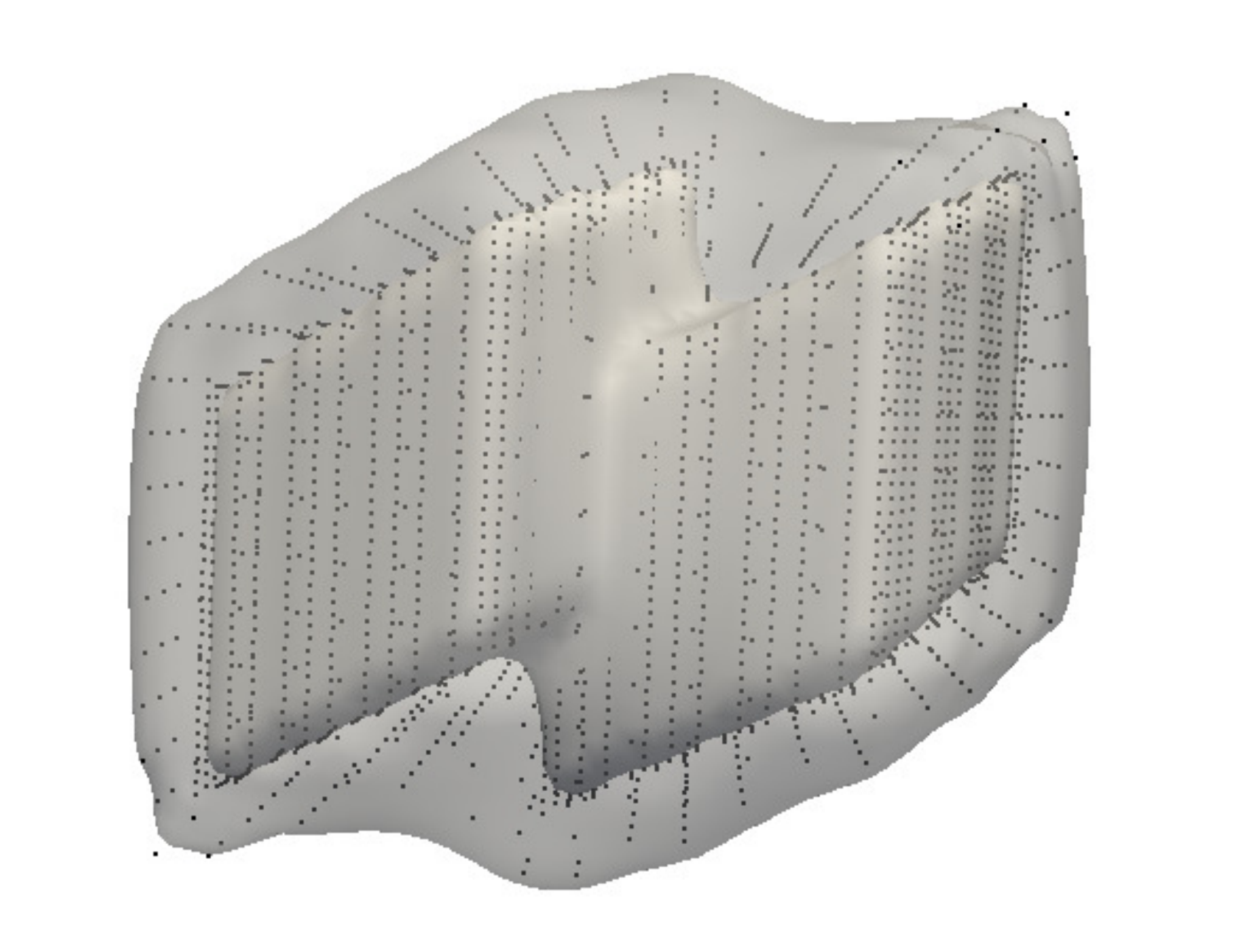,width=8cm}
\caption{{Left: reconstructed surface compared with the data set (black dots). Right: reconstructed surface compared with the grid (black dots) and the zero-level set of the initial condition.}\label{cubigrid}}
\end{figure}

As a third test, we consider a data set $S$ made of 4020 points uniformly chosen on a shape made by two intersecting cubes in $[-2,2]^3$, with axes not aligned with the grid. We consider a grid $\G$ containing  the data set and  $1076$ additional nodes uniformly  chosen in a narrow band, computed according to \eqref{n_band} with $\delta_S=0.01$ from  a full grid of $80\times 80\times 80$ nodes. We apply \eqref{Scheme3d} on the grid $\G$ with time step $\D t=0.01$. The algorithm is stopped after 100 iteration, and the final solution is shown in the left plot of Fig. \ref{cubisol}. Data set, grid points and initial condition are shown in Fig. \ref{cubigrid}.

\paragraph{Teapot}

\begin{figure}
\centering
\epsfig{figure=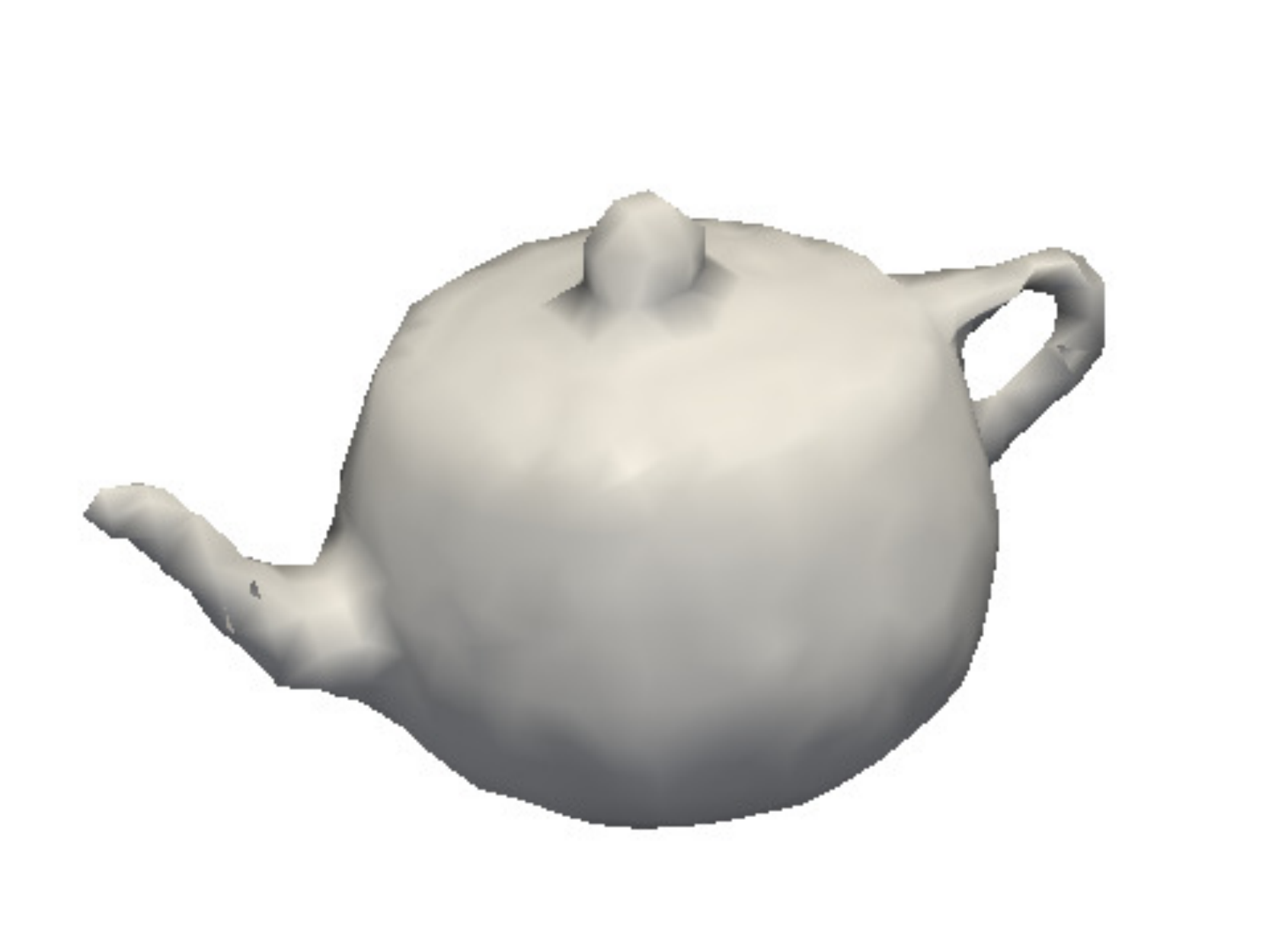,width=8cm}\epsfig{figure=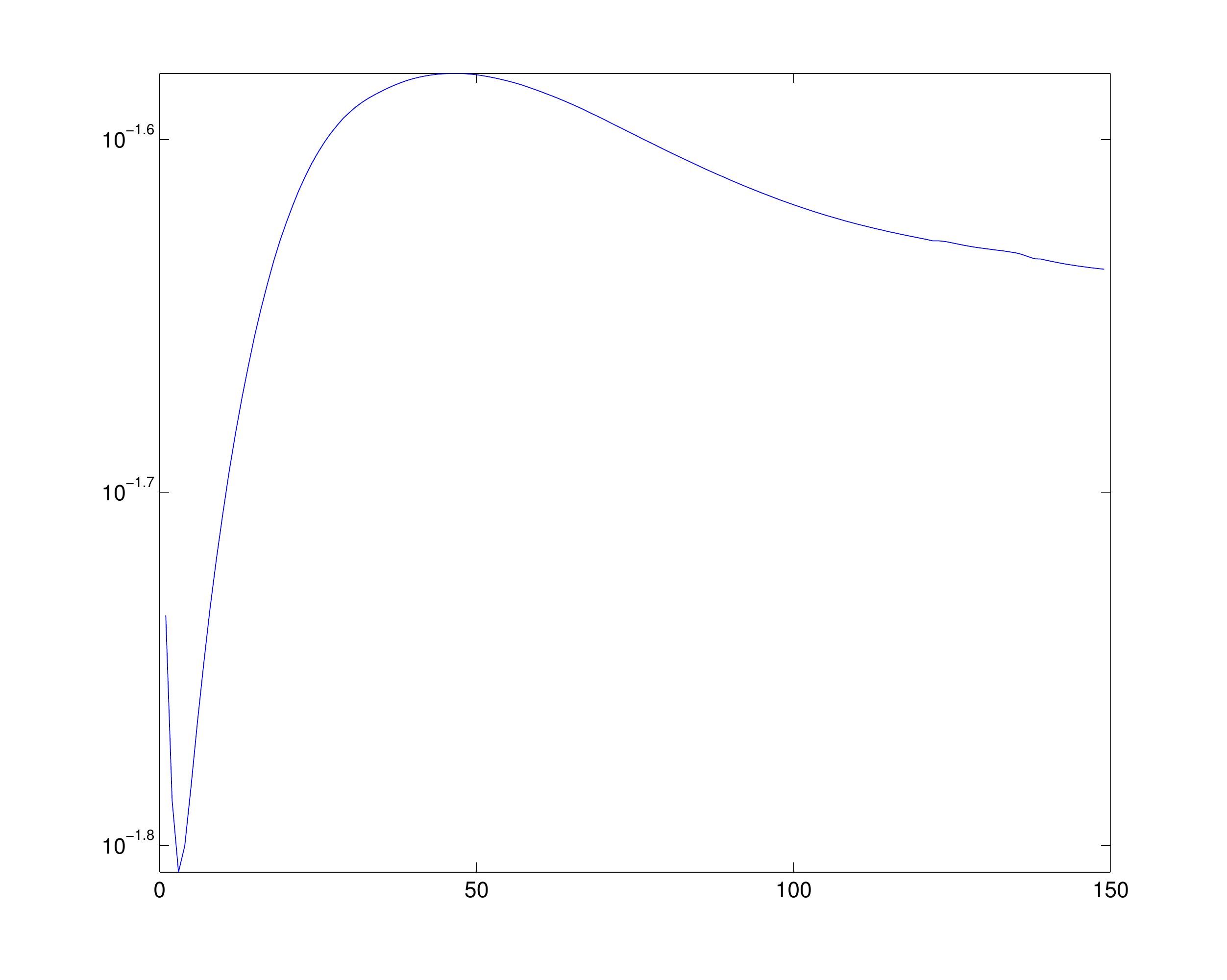,width=8cm}
\caption{{Left: reconstructed surface for the teapot. Right: relative update between two successive iterations}
\label{Test7sol}}
\epsfig{figure=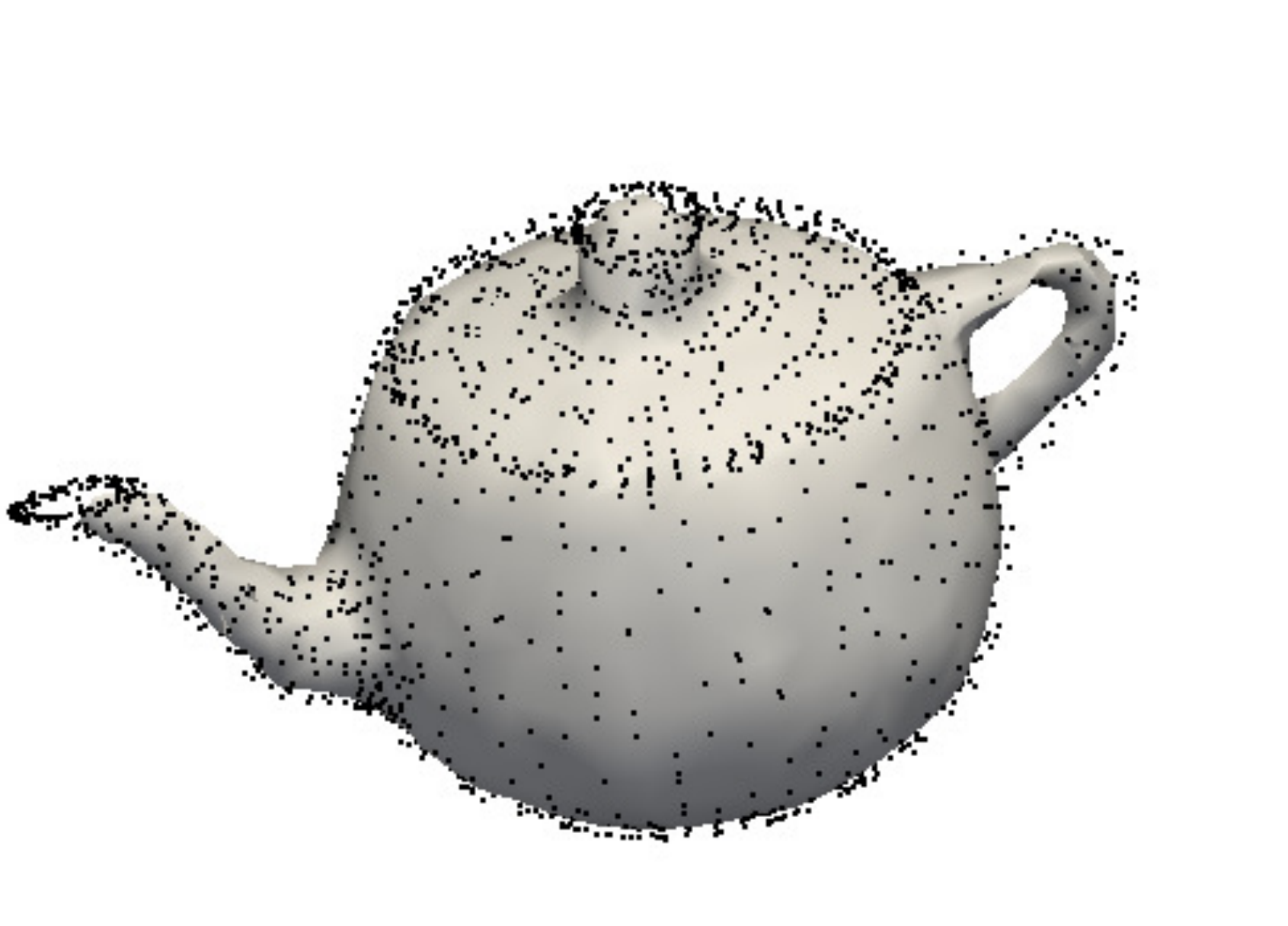,width=8cm}\epsfig{figure=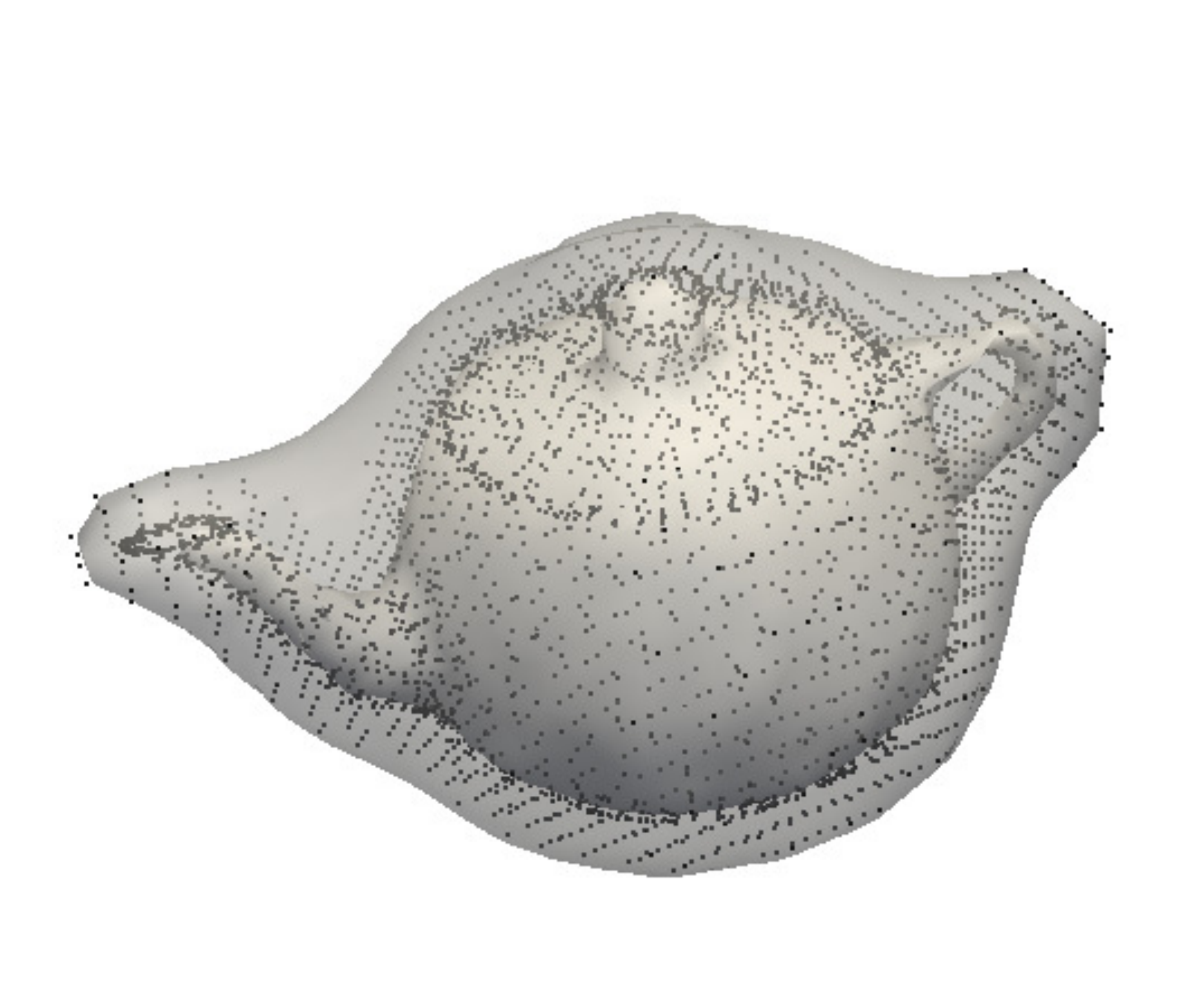,width=8cm}
\caption{{Left : reconstructed surface compared with the data set (black dots). Right: reconstructed surface compared with the grid (black dots) and the zero-level set of the initial condition.}\label{Test7grid}}
\end{figure}

Finally, we test our algorithm on the more complex benchmark of a teapot-shaped surface (whose data set has been taken from {\url{http://www.itl.nist.gov/iad/vug/sharp/benchmark/3DInterestPoint/}}) in the cube $[-0.8,0.8]^3$. 
We  consider a data set $S$ made of 2602 points (about $25\%$ of the original dataset) and a  grid $\G$ containing the dataset and  $3109$ additional nodes chosen on a narrow band of the data set, computed with $\delta_S=0.1$ from a full grid of $50\times50\times50$ nodes.  We apply \eqref{Scheme3d} on the grid $\G$ with time step $\D t=0.001$. The algorithm is stopped after 150 iterations, showing the final solution and the convergence history in the two plots of Fig. \ref{Test7sol}, while the data set and the grid are shown in Fig. \ref{Test7grid}.

\end{document}